\documentclass[seceqn,dvips]{arxbj}
\usepackage{upgreek}
\usepackage{graphicx}

\aid{0}
\volume{16}
\issue{3}
\pubyear{2010}
\firstpage{730}
\lastpage{758}
\doi{10.3150/09-BEJ231}

\makeatletter

\newtheorem{theorem}{Theorem}
\newtheorem{lem}{Lemma}
\newremark{remark}{Remark}
\newtheorem{coro}{Corollary}
\newtheorem{prop}{Proposition}
\newremark{example}{Example}

\makeatother

\begin{document}
\begin{frontmatter}

\title{Varying-coefficient functional linear regression}
\runtitle{Varying-coefficient functional linear regression}

\begin{aug}
\author[1]{\fnms{Yichao} \snm{Wu}\corref{}\thanksref{1}\ead[label=e1]{wu@stat.ncsu.edu}},
\author[2]{\fnms{Jianqing} \snm{Fan}\thanksref{2}\ead[label=e2]{jqfan@princeton.edu}}
\and
\author[3]{\fnms{Hans-Georg} \snm{M\"{u}ller}\thanksref{3}\ead[label=e3]{mueller@wald.ucdavis.edu}}
\runauthor{Y. Wu, J. Fan and H.-G. M\"{u}ller}
\address[1]{Department of Statistics, North Carolina State University,
Raleigh, NC 27695, USA.\\ \printead{e1}}
\address[2]{Department of Statistics, Princeton University, Princeton,
NJ 08544, USA.\\ \printead{e2}}
\address[3]{Department of Statistics, University of California--Davis,
Davis, CA 95616, USA.\\ \printead{e3}}
\end{aug}

\received{\smonth{8} \syear{2008}}
\revised{\smonth{4} \syear{2009}}

%
\begin{abstract}
Functional linear regression analysis aims to model regression
relations which include a functional predictor. The analog of the
regression parameter vector or matrix in conventional
multivariate or multiple-response linear regression models is a regression
parameter function in one or two arguments. If, in addition, one has
scalar predictors, as is often the case in applications to
longitudinal studies, the question arises how to incorporate these into
a functional regression model. We study a varying-coefficient
approach where the scalar covariates are modeled as additional
arguments of the regression parameter function. This extension of
the functional linear regression model is analogous to the extension
of conventional linear regression models to varying-coefficient
models and shares its advantages, such as increased flexibility;
however, the details of this extension are more challenging in the
functional case. Our methodology combines smoothing methods with
regularization by truncation at a finite number of functional
principal components. A practical version is developed and
is shown to perform better than functional linear regression for
longitudinal data. We investigate the asymptotic properties of
varying-coefficient functional linear regression and establish
consistency properties.
\end{abstract}

%

\end{frontmatter}

\section{Introduction}

Functional linear regression analysis is an extension of ordinary
regression to the case where predictors include random functions and
responses are scalars or functions. This methodology has recently
attracted increasing interest due to its inherent applicability in
longitudinal data analysis and other areas of modern data analysis.
For an excellent introduction, see Ramsay and Silverman (\citeyear{RS2005}).
Assuming that predictor process $X$ possesses a square-integrable
trajectory (i.e., $X\in L^2(\mathcal{S})$, where
$\mathcal{S}\subset\mathbb{R}$), commonly considered functional linear
regression models include
%
\begin{equation}\label{scalarflr}
E(Y|X)=\mu_Y+\int_{\mathcal{S}}\beta(s)\bigl(X(s)-\mu_X(s)\bigr)\,\mathrm{d}s ,
\end{equation}
with a scalar response $Y\in \mathbb{R}$, and\vspace{-2pt}
%
\begin{equation}\label{fcnflr}
E(Y(t)|X)=\mu_Y(t)+\int_{\mathcal{S}}\beta(s,t)\bigl(X(s)-\mu_X(s)\bigr)\,\mathrm{d}s,
\end{equation}
with a functional response $Y\in L^2(\mathcal{T})$ and
$\mathcal{T}$ being a subset of the real line $\mathbb{R}$, where
$\mu_{X}(s)=E(X(s))$, $s\in\mathcal{S}$ and $\mu_{Y}(t)=E(Y(t))$,
$t\in\mathcal{T}$ (\citet{ramsaydalzell1991}). In analogy to the
classical regression case, estimating equations for the regression
function are based on minimizing the deviation
\[
\beta^*(s, t)=\operatorname{argmin}\limits_{\beta\in
L_2(\mathcal{S}\times\mathcal{T})} E\biggl\{\int_{\mathcal{T}}
\biggl(Y(t)-\mu_Y(t)-\int_{\mathcal{S}}\beta(s,
t)[X(s)-\mu_X(s)]\,\mathrm{d}s\biggr)^2\,\mathrm{d}t\biggr\},
\]
and analogously for (\ref{scalarflr}). To provide a regularized
estimator, one approach is to expand $\beta(\cdot, \cdot)$ in terms
of the eigenfunctions of the covariance functions
of $X$ and $Y$, and to use an appropriately chosen finite number of
the resulting functional principal component (FPC) scores of $X$ as
predictors; see, for example, \citet{MR1389877}, \citeauthor{RS2002} (\citeyear{RS2002,RS2005}), \citet{MR848110},
\citet{ramsaydalzell1991}, \citet{MR1094283}, \citet{MR1789811}, \citet{CFF2002}, \citet{CFS2003SS},
\citet{HH2007}, \citet{CH2006}, \citet{Cardot2007} and many
others.

Advances in modern technology enable us to collect massive amounts
of data at fairly low cost. In such settings, one may observe scalar
covariates, in addition to functional predictor and response
trajectories. For example, when predicting a response such as blood
pressure from functional data, one may wish to utilize functional
covariates, such as body mass index, and also
additional non-functional covariates $Z$, such as the age of a subject.
It is often realistic to
expect
the regression relation to change as an additional covariate such as
age varies. To broaden the applicability of functional linear
regression models, we propose to generalize models (\ref{scalarflr})
and (\ref{fcnflr}) by allowing the slope function to depend on some
additional scalar covariates $Z$. Previous work on
varying-coefficient functional regression models, assuming the case
of a scalar response and of continuously observed predictor
processes, is due to \citet{cardotsarda2005} and recent investigations
of the varying-coefficient approach include \citet{FHL2007} and \citet{ZLJOD2007}.

For ease of presentation, we consider
the case of a one-dimensional covariate $Z\in\mathcal{Z}\subset
\mathbb{R}$,
extending (\ref{scalarflr}) and (\ref{fcnflr})
to the varying-coefficient functional linear regression
models\vspace*{-2pt}
%
\begin{eqnarray}\label{scalarflrvc}
E(Y|X, Z)=\mu_{Y|Z}+\int_{\mathcal{S}}\beta(Z,
s)\bigl(X(s)-\mu_{X|Z}(s)\bigr)\,\mathrm{d}s
\end{eqnarray}\vspace*{-2pt}
and\vspace*{-2pt}
%
\begin{eqnarray}\label{fcnflrvc}
E(Y(t)|X, Z)=\mu_{Y|Z}(t)+\int_{\mathcal{S}}\beta(Z, s,
t)\bigl(X(s)-\mu_{X|Z}(s)\bigr)\,\mathrm{d}s\vspace{-2pt}
\end{eqnarray}
for scalar and functional responses, respectively, with
corresponding characterizations for the regression parameter
functions\vspace{-2pt}
\begin{eqnarray*}
\beta^*(z, s)&=&\mathop{\operatorname{argmin}}_{\beta(z, \cdot
)\in
L_2(\mathcal{S})} E\biggl\{ \biggl(Y-\mu_{Y|Z}
-\int_{\mathcal{S}} \beta(Z,
s) [X(s)-\mu_{X|Z}(s)]\,\mathrm{d}s\biggr)^2\Bigm| Z=z\biggr\},
\\[-2pt]
\beta^*(z, s, t)&=&\mathop{\operatorname{argmin}}_{\beta(z, \cdot
, \cdot)\in
L_2(\mathcal{S}\times\mathcal{T})} E\biggl\{\int_{\mathcal{T}}
\biggl(Y(t)-\mu_{Y|Z}(t)
\\[-2pt]
&&{}\hspace{94pt}-\int_{\mathcal{S}}\beta(Z, s,
t)[X(s)-\mu_{X|Z}(s)]\,\mathrm{d}s\biggr)^2\,\mathrm{d}t\Bigm|{Z=z}\biggr\}.
\end{eqnarray*}
Here, $\mu_{X|Z}(s)$ and $\mu_{Y|Z}(t)$ denote the conditional mean
function of $X$ and $Y$, given $Z$.

Intuitively, after observing a sample of $n$ observations, $\{
X_i, Y_i, Z_i\}_{i=1}^{n}$, the estimation of the varying
slope functions can be achieved using kernel methods, as follows:
\begin{eqnarray*}
\tilde\beta^*(z, s)&=&\operatorname{argmin}\sum
_{i=1}^{n} K_b(Z_i-z)
\biggl[Y_i-\mu_{Y|Z_i} -\int_{\mathcal{S}} \beta(Z_i, s)
[X_i(s)-\mu_{X|Z_i}(s)]\,\mathrm{d}s\biggr]^2
\end{eqnarray*}
and
\begin{eqnarray*}
\tilde\beta^*(z, s, t)&=&\operatorname{argmin} \sum_{i=1}^{n}
K_b({Z_i-z})
\\
&&{}\hspace{43pt}\times\int_{\mathcal{T}}
\biggl[Y_i(t)-\mu_{Y|Z_i}(t)-\int_{\mathcal{S}}\beta(Z_i, s,
t)[X_i(s)-\mu_{X|Z_i}(s)]\,\mathrm{d}s\biggr]^2\,\mathrm{d}t
\end{eqnarray*}
for (\ref{scalarflrvc}) and (\ref{fcnflrvc}), respectively, where
$K_b(z)=K(z/b)/b$ for a kernel function $K(\cdot)$ and a bandwidth
$b>0$. The necessary regularization of the slope function is
conveniently achieved by truncating the Karhunen--Lo\`{e}ve expansion
of the covariance function for the predictor process (and the
response process, if applicable). To avoid difficult technical
issues and enable straightforward and rapid implementation, it is
expedient to adopt the two-step estimation scheme proposed and
extensively studied by \citet{FZ2000}.

To this end, we first bin our observations according to the values
taken by the additional covariate $Z$ into a partition of
$\mathcal{Z}$. For each bin, we obtain the sample covariance
functions based on the observations within this bin. Assuming that
the covariance functions of the predictor and response processes are
continuous in $z$ guarantees that these sample covariance functions
converge to the corresponding true covariance functions evaluated at
the bin centers as bin width goes to zero and sample size increases.
This allows us to estimate the slope function at each bin center
consistently, using the technique studied in \citet{YMW2005annal},
providing initial raw estimates. Next, local linear smoothing
(\citet{FG1996}) is applied to improve estimation efficiency,
providing our final estimator of the slope function for any
$z\in\mathcal{Z}$.

The remainder of the paper is organized as follows. In Section
\ref{sec:estimation}, we introduce basic notation and present our
estimation scheme. Asymptotic consistency properties are reported in
Section \ref{sec:asymptotics}.
Finite-sample implementation issues are discussed in Section
\ref{sec:fsi}, results of simulation studies in Section
\ref{sec:simulation} and real data applications in Section
\ref{sec:application}, with conclusions in Section 7. Technical
proofs and auxiliary results are given in the
\hyperref[app]{Appendix}.

\vspace*{-3pt}
\section{Varying coefficient functional linear regression for
sparse and irregular data} \label{sec:estimation}
\vspace*{-3pt}

To facilitate
the presentation, we focus on the case of a functional response,
which remains largely unexplored. The case with a scalar response
can be handled similarly. We also emphasize the case of sparse and
irregularly observed data with errors, due to its relevance in
longitudinal studies. The motivation of the varying-coefficient
functional regression models (\ref{scalarflrvc}) and
(\ref{fcnflrvc}) is to borrow strength across subjects, while
adequately reflecting the effects of the
additional covariate. We impose the following smoothness conditions:
\begin{enumerate}[{[A0]}]
\item[{[A0]}] The conditional mean and covariance functions of the
predictor and response processes depend on $Z$ and are continuous in
$Z$, that is, $\mu_{X,z}(s)=E(X(s)|Z=z)$,
$\mu_{Y,z}(t)=E(Y(t)|Z=z)$, $G_{X,z}(s_1, s_2)=\operatorname{cov}(X(s_1),
X(s_2)|Z=z)$, $G_{Y,z}(t_1, t_2)=\operatorname{cov}(Y(t_1),
Y(t_2)|Z=z)$ and
$C_{XY,z}(s, t)=\operatorname{cov}(X(s), Y(t)|Z=z)$
are continuous in $z$ and their respective arguments, and
have continuous second order partial derivatives with respect to~$z$.
\end{enumerate}

Note that [A0] implies that the conditional mean and covariance
functions of predictor and response processes do not change radically
in a small neighborhood of $Z=z$. This facilitates the estimation of
$\beta(z, s, t)$, using the two-step
estimation scheme proposed by \citet{FZ2000}. While, there,
the additional covariate $Z$ is assumed to take values on a grid,
in our case, $Z$ is more generally assumed to be continuously
distributed. In this case, we assume that the additional variable $Z$
has a compact domain $\mathcal{Z}$ and its density $f_Z(z)$ is
continuous and bounded away from both zero and infinity.

\begin{enumerate}[{[A1]}]
\item[{[A1]}] $\mathcal{Z}$ is compact, $f_Z(z)\in C^0$, $\underline
{f_{Z}} =\inf_{z\in\mathcal{Z}} f_Z(z)>0$ and $\bar{f_Z}= \sup
_{z\in\mathcal{Z}} f_Z(z)<\infty$.
\end{enumerate}

\subsection{Representing predictor and response functions via
functional principal components for sparse and irregular data}
Suppose that we have observations on $n$ subjects. For each subject $i$,
conditional on $Z_i=z_i$, the square-integrable predictor trajectory
$X_i$ and response trajectory $Y_i$ are unobservable realizations of
the smooth random processes $(X, Y|Z=z_i)$, with unknown mean and
covariance functions (condition [A0]). The arguments of $X(\cdot)$
and $Y(\cdot)$ are usually referred to as time. Without loss of
generality, their domains $\mathcal{S}$ and $\mathcal{T}$ are
assumed to be finite and closed intervals. Adopting the general
framework of functional data analysis, we assume, for each $z$, that
there exist orthogonal expansions of the covariance functions
$G_{X,z}(\cdot, \cdot)$ (resp.~$G_{Y,z}(\cdot, \cdot)$) in the $L_2$
sense via the eigenfunctions $\psi_{z,m}$ (resp.~$\phi_{z,k}$), with
non-increasing eigenvalues $\rho_{z,m}$ (resp.~$\lambda_{z,k}$), that is,
$G_{X,z}(s_1, s_2)=\sum_{m=1}^\infty
\rho_{z,m}\psi_{z,m}(s_1)\psi_{z,m}(s_2)$, $G_{Y,z}(t_1,
t_2)=\sum_{k=1}^\infty\lambda_{z,k}\phi_{z,k}(t_1)\phi_{z,k}(t_2)$.

Instead of observing the full predictor trajectory $X_i$ and
response trajectory $Y_i$, typical longitudinal data consist of
noisy observations that are made at sparse and irregularly spaced
locations or time points, providing sparse measurements of predictor
and response trajectories that are contaminated with additional
measurement errors (\citet{staniswalislee1998}, \citet{ricewu2001},
\citeauthor{YMW2005jasa}  (\citeyear{YMW2005jasa,YMW2005annal})). To adequately reflect the situation of
sparse, irregular and possibly subject-specific time points
underlying these measurements, we assume that a random number $L_i$
(resp.~$N_i$) of measurements for $X_i$ (resp.~$Y_i$) is made, at
times denoted by $S_{i1}, S_{i2}, \dots, S_{iL_i}$ (resp.~$T_{i1},
T_{i2}, \dots, T_{iN_i}$). Independent of any other random
variables, the numbers of points sampled from each trajectory
correspond to random variables $L_i$ and $N_i$ that are assumed to
be i.i.d.~as $L$ and $N$ (which may be correlated), respectively.
For $1\leq i\leq n$, $1\leq l\leq L_i$, $1\leq j\leq N_i$, let
$U_{il}$ (resp.~$V_{ij}$) be the observation of the random
trajectory $X_i$ (resp.~$Y_i$) made at a random time $S_{il}$
(resp.~$T_{ij}$), contaminated with measurement errors $\varepsilon_{il}$
(resp.~$\epsilon_{ij}$).
Here, the random measurement errors $\varepsilon_{il}$ and
$\epsilon_{ij}$ are assumed to be i.i.d., with mean zero and variances
$\sigma_X^2$ and $\sigma_Y^2$, respectively. They are independent of
all other random variables. The following two assumptions are
made.

\begin{enumerate}[{[A2]}]
\item[{[A2]}] For each subject $i$,
$L_i\stackrel{\mathrm{i.i.d.}}{\sim}L$ (resp.~$N_i\stackrel{\mathrm{i.i.d.}}{\sim}N$) for
a positive discrete-valued random variable with $EL<\infty$
(resp.~$EN<\infty$) and $P(L>1)>0$ (resp.~$P(N>1)>0$).

\item[{[A3]}] For each subject $i$, observations
on $X_i$ (resp.~$Y_i$) are independent of $L_i$ (resp.~$N_i$), that is,
$\{(S_{il},
U_{il}\dvt l\in\mathcal{L}_i)\}$ is independent of $L_i$ for any
$\mathcal{L}_i\subset\{1, 2, \ldots, L_i\}$ (resp.~$\{(T_{ij},
V_{ij})\dvt j\in\mathcal{N}_i\}$ is independent of $N_i$ for any
$\mathcal{N}_i\subset\{1, 2, \ldots, N_i\}$).
\end{enumerate}

It it surprising that under these ``longitudinal assumptions'', where
the number of observations per subject is fixed and does not
increase with sample size, one can nevertheless obtain asymptotic
consistency results for the regression relation. This phenomenon was
observed in \citet{YMW2005annal} and is due to the fact that,
according to (\ref{betarepresentation}), the target regression
function depends only on localized eigenfunctions, localized
eigenvalues and cross-covariances of localized functional principal
components. However, even though localized, these eigenfunctions and
moments can be estimated from pooled data and do not require the
fitting of individual trajectories. Even for the case of fitted
trajectories, conditional approaches have been implemented
successfully, even allowing reasonable derivative
estimates to be obtained from very sparse data (\citet{Liumueller2009}).

Conditional on $Z_i=z$, the FPC scores of $X_i$ and $Y_i$ are
$\zeta_{z, im}=\int[X_i(s)-\mu_{X,z}(s)]\psi_{z,m}(s) \,\mathrm{d}s$ and
$\xi_{z,ik} =\int[Y_i(s)-\mu_{Y,z}(s)]\phi_{z,k}(s) \,\mathrm{d}s,$
respectively. For all $z$, these FPC scores $\zeta_{z, im}$ satisfy
$E\zeta_{z, im}=0$, $\operatorname{corr}(\zeta_{z, im_1}, \zeta_{z,
im_2})=0$ for
any $m_1\ne m_2$ and $\operatorname{var}(\zeta_{z, im})=\rho_{z,
m}$; analogous
results hold
for $\xi_{z,ik}$. With this notation, using the Karhunen--Lo\`{e}ve
expansion as in \citet{YMW2005annal}, conditional on $Z_i$, the
available measurements of the $i$th predictor and response
trajectories can be represented as
\begin{eqnarray}\label{xsparseassum}
U_{il}&=&X_i(S_{il})+\varepsilon_{il}\nonumber
\\[-8pt]\\[-8pt]
&=&\mu_{X,Z_i}(S_{il})+\sum_{m=1}^\infty\zeta_{Z_i,
im}\psi_{Z_i,m}(S_{il})+\varepsilon_{il},   \qquad      1\leq
l\leq L_i,\nonumber
\\\label{ysparseassum}
V_{ij}&=&Y_i(T_{ij})+\epsilon_{ij}\nonumber
\\[-8pt]\\[-8pt]
&=&\mu_{Y,Z_i}(T_{ij})+\sum_{k=1}^\infty\xi_{Z_i,
ik}\phi_{Z_i,k}(T_{ij})+\epsilon_{ij},     \qquad    1\leq j\leq\nonumber
N_i.
\end{eqnarray}
%

\subsection{Estimation of the slope function}
 For estimation of the slope function, one standard approach is
to expand it in terms of an orthonormal functional basis and to
estimate the coefficients of this expansion to estimate the slope
function in the non-varying model (\ref{fcnflr})
(\citet{YMW2005annal}). As a result of the non-increasing property of
the eigenvalues of the covariance functions, such expansions of the
slope function are often efficient and require only a few components
for a good approximation. Truncation at a finite number of terms
provides the necessary regularization. Departing from
\citet{YMW2005annal}, we assume here that an additional covariate $Z$
plays an important role and must be incorporated into the model,
motivating (\ref{fcnflrvc}). To make this model as flexible as
possible, the conditional mean and covariance functions of the
predictor and response processes are allowed to change smoothly with
the value of the covariate $Z$ (Assumption [A0]), which facilitates
implementation and analysis of the two-step estimation scheme, as in
\citet{FZ2000}.

Efficient implementation of the two-step estimation scheme begins by
binning subjects according to the levels of the additional covariate
$Z_i$, $i=1,2,\dots, n$. For ease of presentation, we use
bins of equal width, although, in practice, non-equidistant bins can
occasionally be advantageous. Denoting the bin centers by $z^{(p)},
p=1, 2, \dots, P$, and the bin width by $h$, the $p$th bin is
$[z^{(p)}-\frac{h}{2}, z^{(p)}+\frac{h}{2})$ with
$h=\frac{\mid\mathcal{Z}\mid}{P}$, where $|\mathcal{Z}|$ denotes the
size of the domain of $Z$, $z^{(1)}-h/2\equiv\inf\{z\dvt
z\in\mathcal{Z} \}$ and $z^{(P)}+h/2\equiv\sup\{z\dvt
z\in\mathcal{Z}\}$ (note that the last bin is $[z^{(P)}-h/2,
z^{(P)}+h/2]$). Let $\mathcal{N}_{z, h}=\{i\dvt
Z_i\in[z-\frac{h}{2}, z+\frac{h}{2})\}$ be the index set of
those subjects falling into bin $[z-\frac{h}{2}, z+\frac{h}{2})$ and
$n_{z, h}=\# \mathcal{N}_{z, h}$ the number of those subjects.

\subsubsection{Raw estimates}

For each bin $[z^{(p)}-\frac{h}{2},
z^{(p)}+\frac{h}{2})$, we use the
\citet{YMW2005jasa} method to obtain
our raw estimates
$\tilde\mu_{X,z^{(p)}}(\cdot)$ and
$\tilde\mu_{Y,z^{(p)}}(\cdot)$ of the
conditional mean trajectories and the raw
slope function estimate
$\tilde\beta(z^{(p)}, s, t)$. The
corresponding raw estimates of
$\sigma_X^2$ and $\sigma_Y^2$ are
denoted by $\tilde\sigma_{X,
z^{(p)}}^2$ and $\tilde\sigma_{Y,
z^{(p)}}^2$ for $p=1, 2, \dots, P$.
For each $1\leq p \leq P$, the local linear scatterplot smoother of
$\tilde\mu_{X, z^{(p)}}(s)$ is defined through minimizing
\begin{eqnarray}
\sum_{i\in\mathcal{N}_{z^{(p)}, h}} \sum_{j=1}^{N_i}
\kappa_1\biggl(\frac{S_{ij}-s}{b_{X,z^{(p)}}}\biggr)\bigl(U_{ij}-
d_{0}-d_{1}(S_{ij}-s)\bigr)^2\nonumber
\end{eqnarray}
with respect to $d_0$ and $d_1$, and setting $\tilde\mu_{X,
z^{(p)}}(s)$ to be the minimizer $\hat d_0$, where $\kappa_1(\cdot)$
is a kernel function and $b_{X,z^{(p)}}$ is the smoothing bandwidth,
the choice of which will be discussed in Section~\ref{sec:fsi}. We
define a similar local linear scatterplot smoother of $\tilde\mu_{Y,
z^{(p)}}(t)$. According to Lemma~\ref{ymwjasathm11dmean} in the
\hyperref[app]{Appendix}, raw
estimates $\tilde\mu_{X, z^{(p)}}(s)$ and $\tilde\mu_{Y,
z^{(p)}}(t)$ are consistent uniformly for $z^{(p)}$, $p=1, 2,
\dots, P$, for appropriate bandwidths $b_{X,z^{(p)}}$ and $b_{Y,z^{(p)}}$.

Extending \citet{YMW2005annal}, the conditional slope function can be
represented as
%
\begin{eqnarray}\label{betarepresentation}
 \beta(z, s,
t)=\sum_{k=1}^\infty\sum_{m=1}^\infty\frac{E\zeta_{z, m}
\xi_{z,k}}{E\zeta_{z,m}^2}\psi_{z, m}(s)\phi_{z,
k}(t)
\end{eqnarray}
for each $z$, where $\psi_{z,m}(\cdot)$ and $\phi_{z,k}(\cdot)$ are
the eigenfunctions of covariance functions $G_{X,z}(\cdot, \cdot)$
and $G_{Y,z}(\cdot, \cdot)$, respectively, and $\zeta_{z, m}$ and
$\xi_{z,k}$ are the functional principal component scores of $X$ and
$Y$, respectively, conditional on $Z=z$.

To obtain raw slope estimates $\tilde\beta(z^{(p)}, s, t)$ for
$p=1, 2, \dots, P$, we first estimate the conditional covariance
functions $ G_{X,z^{(p)}}(s_1, s_2)$, $
G_{Y,z^{(p)}}(t_1, t_2)$ and $ C_{XY,z^{(p)}}(s, t)$ at each
bin center, based on the observations falling into the bin, using
the approach of \citet{YMW2005annal}. From ``raw'' covariances,
$G_{X,i,z^{(p)}}(S_{ij}, S_{ik})=(U_{ij}-\tilde
\mu_{X,z^{(p)}}(S_{ij}))(U_{ik}-\tilde\mu_{X, z^{(p)}}(S_{ik}))$
for $1\leq j, k\leq L_i$, $i\in\mathcal{N}_{z^{(p)},h}$ and $p=1, 2,
\dots, P$, and the locally smoothed conditional covariance $\tilde
G_{X,z^{(p)}}(s_1, s_2)$ is defined as the minimizer $\hat b_0$ of
the local linear problem
\begin{eqnarray*}
&&\min_{b_0, b_{11}, b_{12}}\sum_{i\in\mathcal{N}_{z^{(p)},h}}
\sum_{1\leq j\ne l\leq L_i}
\kappa_2\biggl(\frac{S_{ij}-s_1}{h_{X,z^{(p)}}},
\frac{S_{il}-s_2}{h_{X,z^{(p)}}}\biggr)
\\
&&{}\hspace{101pt}\times [ G_{X, i,z^{(p)}}(S_{ij},
S_{il})-b_0-b_{11}(S_{ij}-s_1)-b_{12}(S_{il}-s_2)]^2,
\end{eqnarray*}
where $\kappa_2(\cdot, \cdot)$ is a bivariate kernel function and
$h_{X, z^{(p)}}$ a smoothing bandwidth. The diagonal ``raw''
covariances $G_{X, i,z^{(p)}}(S_{ij}, S_{ij})$ are removed from the
objective function of the above minimization problem because
$EG_{X, i,z^{(p)}}(S_{ij}, S_{il})\approx\operatorname{cov}(X(S_{ij}),
X(S_{il}))+\delta_{jl}\sigma_X^2$, where $\delta_{jl}=1$ if $j=l$ and
$0$ otherwise. Analogous considerations apply for $\tilde
G_{Y,z^{(p)}}(T_{ij}, T_{il})$. The diagonal ``raw'' covariances
$G_{X, i,z^{(p)}}(S_{ij}, S_{ij})$ and $G_{Y, i,z^{(p)}}(T_{ij},
T_{ij})$
can be smoothed with bandwidths $b_{X, z^{(p)}, V}$ and $b_{Y, z^{(p)},
V}$, respectively, to estimate $V_{X, z^{(p)}}(s)=
G_{X, z^{(p)}}(s, s)+\sigma_X^2$ and $V_{Y, z^{(p)}}(t)= G_{Y, z^{(p)}}(t,
t)+\sigma_Y^2$, respectively. The resulting estimators are denoted by
$\tilde
V_{X, z^{(p)}}(s)$ and $\tilde V_{Y,z^{(p)}}(t)$, respectively, and
the differences $(\tilde V_{X, z^{(p)}}(s)-\tilde G_{X,z^{(p)}}(s,
s))$ (and analogously for $Y$) can be used to obtain estimates
$\tilde\sigma_{X, z^{(p)}}^2$ for $\sigma_X^2$ and $\tilde\sigma_{Y,
z^{(p)}}^2$ for $\sigma_Y^2$, by integration. Furthermore, ``raw''
conditional cross-covariances $C_{i, z^{(p)}}(S_{il},
T_{ij})=(U_{il}-\tilde\mu_{X,z^{(p)}}(S_{il}))(V_{ij}-\tilde
\mu_{Y,z^{(p)}}(T_{ij}))$ are used to estimate $ C_{XY,z^{(p)}}(s,
t),$ by minimizing
\begin{eqnarray*}
&&\sum_{i\in\mathcal{N}_{z^{(p)},h}} \sum_{1\leq l\leq L_i}
\sum_{1\leq j \leq N_i} \kappa_2\biggl(\frac{S_{ij}-s}{h_{1,z^{(p)}}},
\frac{T_{ij}-t}{h_{2,z^{(p)}}}\biggr)
\\
&&{}\hspace{88pt}\times [ C_{i, z^{(p)}}(S_{il},
T_{ij})-b_0-b_{11}(S_{il}-s)-b_{12}(T_{ij}-t)]^2
\end{eqnarray*}
with respect to $b_0$, $b_{11}$ and $b_{12},$ and setting $\tilde
C_{XY,z^{(p)}}(s, t)$ to be the minimizer $\hat b_0$, with smoothing
bandwidths $h_{1,z^{(p)}}$ and $h_{2,z^{(p)}}$.

In (\ref{betarepresentation}), the slope function may be represented
via the eigenvalues and eigenfunctions of the covariance operators.
To obtain the estimates $\tilde\rho_{z^{(p)}, m}$ and
$\tilde\psi_{z^{(p)}, m}(\cdot)$ (resp.~$\tilde\lambda_{z^{(p)},
k}$ and
$\tilde\phi_{z^{(p)}, k}(\cdot)$) of eigenvalue--eigenfunction pairs
$\rho_{z^{(p)}, m}$ and
$\psi_{z^{(p)}, m}(\cdot)$ (resp.~$\lambda_{z^{(p)}, k}$ and
$\phi_{z^{(p)}, k}(\cdot)$), we use conditional functional principal
component analysis (CFPCA) for $\tilde G_{X,z^{(p)}}(\cdot, \cdot)$
(resp.~$\tilde G_{Y,z^{(p)}}(\cdot, \cdot)$), by numerically solving
the conditional eigenequations
%
\begin{eqnarray}\label{Xcoveigenexp}
\int_{\mathcal{S}}\tilde G_{X,z^{(p)}}(s_1, s_2)\tilde\psi_{z^{(p)},
m}(s_1)\,\mathrm{d}s_1&=&\tilde\rho_{z^{(p)}, m}\tilde\psi_{z^{(p)},
m}(s_2),      \qquad       m=1, 2, \ldots,
\\\label{Ycoveigenexp}
\int_{\mathcal{T}}\tilde G_{Y,z^{(p)}}(t_1, t_2)\tilde\phi_{z^{(p)},
k}(t_1)\,\mathrm{d}t_1&=&\tilde\lambda_{z^{(p)}, k}\tilde\phi_{z^{(p)},
k}(t_2),       \qquad         k=1, 2,
\ldots.
\end{eqnarray}
Note that we estimate the conditional mean functions and conditional
covariance functions over dense grids of $\mathcal{S}$ and $\mathcal
{T}$. Numerical integrations like the one on the left-hand side of
(\ref{Xcoveigenexp}) are done over these dense grids using the
trapezoid rule. Note, further, that integrals over individual
trajectories are not needed for the regression focus, in that we use
conditional expectation to estimate principal component scores, as in
(\ref{eq4.1}).

Due to the fact that
\begin{eqnarray*}
C_{XY,z}(s, t)&=&\operatorname{cov}\bigl(X(s),
Y(t)|Z=z\bigr)=\sum_{k=1}^\infty\sum_{m=1}^\infty{E(\zeta_{z, m}
\xi_{z,k})}\psi_{z, m}(s)\phi_{z, k}(t),
\end{eqnarray*}
we then obtain preliminary estimates of $\sigma_{z, mk}=E(\zeta_{z,
m} \xi_{z,k})$ at the bin centers $z^{(p)}$, $p=1, 2, \dots, P$, by
numerical integration,
%
\begin{eqnarray}\label{xycoveigenexpand}
\tilde\sigma_{z^{(p)},mk}=\int_{\mathcal{T}}\int_{\mathcal{S}}
\tilde\psi_{z^{(p)},m}(s) \tilde C_{XY,z^{(p)}}(s, t)
\tilde\phi_{z^{(p)},k}(t)\,\mathrm{d}s\,\mathrm{d}t.
\end{eqnarray}
With (\ref{betarepresentation}), (\ref{Xcoveigenexp}),
(\ref{Ycoveigenexp}) and (\ref{xycoveigenexpand}), the raw
estimates of $\beta(z^{(p)}, s, t)$ are
%
\begin{eqnarray}\label{estbetaraw1}
\tilde\beta\bigl(z^{(p)}, s, t\bigr)=\sum_{k=1}^K\sum_{m=1}^M \frac{\tilde
\sigma_{z^{(p)}, {mk}}}{\tilde\rho_{z^{(p)}, m}}
\tilde\psi_{z^{(p)},m}(s) \tilde\phi_{z^{(p)},k}(t).
\end{eqnarray}
Further details on the ``global'' case can be found in
\citet{YMW2005annal}.

\subsubsection{Refining the raw estimates}

We establish in the \hyperref[app]{Appendix} that the raw estimates
$\tilde\mu_{X,z^{(p)}}(s)$, $\tilde\mu_{Y,z^{(p)}}(t)$ and
$\tilde\beta(z^{(p)},s,t)$ are consistent. As has been demonstrated in
\citet{FZ2000}, there are several reasons to refine such raw
estimates. For example, the raw estimates are generally not smooth
and are based on local observations, hence inefficient. Most
importantly, applications require that the function $\beta(z,s,t)$
is available for any $z\in\mathcal{Z}$.

To refine the raw estimates, the classical approach is smoothing, for
which we adopt the local polynomial smoother. Defining $\mathbf{c}_p=(1,
z^{(p)}-z, \dots, (z^{(p)}-z)^r)^T$, $p=1, 2, \dots, P$, the local
polynomial smoothing weights for estimating the $q$th derivative of
an underlying function are
\begin{eqnarray*}
\omega_{q, r+1}\bigl(z^{(p)}, z, b\bigr)= q! \mathbf{e}^T_{q+1, r+1}(\mathbf
{C}^T\mathbf{W}\mathbf{C}
)^{-1}\mathbf{c}_p K_b\bigl(z^{(p)}-z\bigr),
  \qquad  p=1, 2, \ldots, P,
\end{eqnarray*}
where $\mathbf{C}=(\mathbf{c}_1, \mathbf{c}_2, \dots, \mathbf
{c}_P)^T$, $\mathbf{W}
=\operatorname{diag}(K_b(z^{(1)}-z), K_b(z^{(2)}-z), \dots, K_b(z^{(P)}-z))$ and
$\mathbf{e}
_{q+1, r+1}=(0, \dots, 0, 1, 0,
\dots, 0)^T$ is a unit vector of length $r+1$ with the $(q+1)$th
element\vadjust{\goodbreak} being $1$ (see \citet{FG1996}). Our final estimators are
given by
\begin{eqnarray*}
\hat\mu_{X,z}(s)&=&\sum_{p=1}^P \omega_{0,2}\bigl(z^{(p)}, z,
b\bigr)\tilde\mu_{X,z^{(p)}}(s),
\\
\hat\mu_{Y,z}(t)&=&\sum_{p=1}^P \omega_{0,2}\bigl(z^{(p)}, z,
b\bigr)\tilde\mu_{Y,z^{(p)}}(t),
\\
\hat\beta(z,s,t)&=&\sum_{p=1}^P \omega_{0,2}\bigl(z^{(p)}, z,
b\bigr)\tilde\beta\bigl(z^{(p)}, s , t\bigr).
\end{eqnarray*}
Due to the assumption that the variance of the measurement error
does not depend on the additional covariate, the final estimators
of $\sigma_X^2$ and $\sigma_Y^2$ can be taken as simple averages,
%
\begin{eqnarray}\label{errorvarfinal}
\hat\sigma_X^2=\sum_{p=1}^{P}\tilde\sigma_{X, z^{(p)}}^2/P \quad
\mbox{and}\quad
\hat\sigma_Y^2=\sum_{p=1}^{P}\tilde\sigma_{Y, z^{(p)}}^2/P.
\end{eqnarray}

\begin{remark}
The localization to $Z=z$, as needed for the
proposed varying coefficient model,
coupled with the extreme sparseness assumption \textup{[A2]}, which
adequately reflects longitudinal designs, is not conducive to
obtaining explicit results in terms of convergence rates for the
general case.
However, by suitably modifying our arguments and coupling them with the
rates of convergence provided on page \textup{2891} of
\citet{YMW2005annal}, we can obtain rates if desired. These are the
rates given there, which depend on complex intrinsic properties of the
underlying processes,
provided that the sample size $n$ is everywhere replaced by $nh$, the
sample size for each bin.
\end{remark}

\begin{remark}
In this work, we focus on sparse and irregularly observed longitudinal
data. For the case where entire processes are observed without noise and
are error-free, one can estimate the localized eigenfunctions at rates
of $L^2$-convergence of $(n\tilde h)^{-1/2}$ (see \citet{HallMuellerWang2006}), where $\tilde h$ is the smoothing
bandwidth. For the
moments of the functional principal components, a smoothing step is not
needed. Known results will be adjusted by replacing $n$ with
$nh$ when conditioning on a fixed covariate level $Z=z$; see \citet
{CH2006} and \citet{HH2007}.
\end{remark}

\section{Asymptotic properties}\label{sec:asymptotics}

We establish some key asymptotic consistency properties for the
proposed estimators. Detailed technical conditions and proofs can be
found in the \hyperref[app]{Appendix}.

The observed data set is denoted by
$\{Z_i, (S_{il}, U_{il})_{l=1}^{L_i}, (T_{ij},
V_{ij})_{j=1}^{N_i}\dvt i=1, 2, \dots, n\}$. We\vspace{1pt} assume that it comes from
(\ref{fcnflr}) and satisfies [A0], [A1], [A2] and [A3].

For $\tilde n \propto\sqrt{n}$, define the event
%
\begin{equation}\label{aprobdefn}
E_n=\{\min n_{z^{(p)}, h} > \tilde n\},
\end{equation}
where $n_{z^{(p)}, h}$ is the number of observations in the $p$th bin
and $\tilde n \propto\sqrt{n}$ means that there exist $c_0$ and $C_0$
such that $0< c_0 \leq\tilde n/\sqrt{n}\leq C_0 <\infty$. It is
shown in Proposition \ref{aprobprop} in the \hyperref[app]{Appendix} that
$P(E_n)\rightarrow1$ as $n\rightarrow\infty$ for $P\propto n^{1/8}$,
as specified by condition (xi).

The global consistency of the final mean and slope function estimates
follows from the following theorem.
\begin{theorem}[(Consistency of time-varying functional regression)]\label
{theorem:finalbetaconv} Under conditions \textup{[A0]}, \textup{[A1]}, \textup{[A2]} and \textup{[A3]} in
Section~\ref{sec:estimation} and conditions \textup{[A4]},
\textup{[A5]} and \textup{(i)--(xi)} in the \hyperref[app]{Appendix}, on the
event $E_n$ with $P(E_n)\rightarrow1$ as $n\rightarrow\infty$, we have
\begin{eqnarray*}
\int_{\mathcal{Z}}\int_{\mathcal{R}}
\bigl(\hat\mu_{W,z}(r)-\mu_{W,z}(r)\bigr)^2\,\mathrm{d}r\,\mathrm{d}z \stackrel{P}{\rightarrow} 0
\qquad \mbox{for }
W=X, \mathcal{R}=\mathcal{S} \mbox{ and } W=Y, \mathcal{R}=\mathcal
{T},
\end{eqnarray*}
and
\begin{eqnarray*}
\int_{\mathcal{Z}}\int_{\mathcal{T}}\int_{\mathcal{S}}
\bigl(\hat\beta(z,s,t)-\beta(z,s,t)\bigr)^2\, \mathrm{d}s\,\mathrm{d}t\,\mathrm{d}z \stackrel{P}{\rightarrow
} 0.
\end{eqnarray*}
\end{theorem}

To study prediction through time-varying functional regression,
consider a new predictor process $X^*$ with associated covariate
$Z^*$. The corresponding conditional expected response process $Y^*$
and its prediction $\hat Y^*$ are given by
%
\begin{eqnarray}\label{theorem2eq1}
Y^*(t)&=&E(Y(t)|X^*, Z^*)\nonumber
\\[-8pt]\\[-8pt]
&=&\mu_{Y, Z^*}(t)+\int_{\mathcal{S}}
\beta(Z^*, s,
t)\bigl(X^*(s)-\mu_{X, Z^*}(s)\bigr)\,\mathrm{d}s,\nonumber
\\\label{theorem2eq2}
\hat Y^*(t)&=&\hat\mu_{Y, Z^*}(t)+\int_{\mathcal{S}}\hat\beta(Z^*,
s, t)\bigl(X^*(s)-\hat\mu_{X, Z^*}(s)\bigr)\,\mathrm{d}s.
\end{eqnarray}

\begin{theorem}[(Consistency of prediction)]\label{theorem4} For a new
predictor process $X^*$ with associated covariate $Z^*$, it holds under
the conditions of Theorem~\ref{theorem:finalbetaconv}
that $\int_{\mathcal{T}}(Y^*(t)-\hat Y^*(t))^2\,\mathrm{d}t\stackrel
{P}{\rightarrow}0$, where
$Y^*(t)$ and $\hat Y^*(t)$ are given by (\ref{theorem2eq1}) and (\ref{theorem2eq2}).
\end{theorem}

\section{Finite-sample implementation}\label{sec:fsi}

For the finite-sample case, several smoothing parameters need to be
chosen. Following
\citet{YMW2005jasa}, the leave-one-curve-out cross-validation method
can be used to select smoothing parameters $b_{X, z^{(p)}}$, $b_{Y,
z^{(p)}}$, $b_{X, z^{(p)},V}$, $b_{Y, z^{(p)},V}$, $h_{X, z^{(p)}}$,
$h_{Y, z^{(p)}}$, $h_{1, z^{(p)}}$ and $h_{2, z^{(p)}}$,
individually for each bin. Further required choices concern the bin
width $h$, the smoothing bandwidth~ $b$ and the numbers $M$ and $K$
of included expansion terms in (\ref{estbetaraw1}). The method of
cross-validation could also be used for these additional choices,
but this incurs a heavy computational load. A fast alternative is
a pseudo-Akaike information criterion (AIC) (or pseudo-Bayesian
information criterion (BIC)).
\begin{enumerate}[{[1]}]
\item[{[1]}] Choose the number of terms in the truncated double
summation representation $\tilde\beta(z^{(p)}, s, t)$ for $M(n)$
and $K(n)$, using AIC or BIC, as in \citet{YMW2005annal}.

\item[{[2]}] For each bin width $h$, choose the best smoothing
bandwidth $b^*(h)$ by minimizing
AIC or BIC.

\item[{[3]}] Choose the bin width $h^*$ which minimizes AIC or BIC,
while, for each $h$
investigated, we use $b^*(h)$ for $b$.
\end{enumerate}

For [1], we will choose $M$ and $K$ simultaneously for all bins,
minimizing the conditional penalized pseudo-deviance given by
\begin{eqnarray*}
C(K)&=&\sum_{p=1}^{P}\sum_{i\in\mathcal{N}_p}\biggl\{\frac
{1}{\tilde\sigma_{Y,
z^{(p)}}^2}\tilde{\bolds\epsilon}_i^T
\tilde{\bolds\epsilon}_i+{N_i}\log
(2\uppi)+{N_i}\log\tilde\sigma_{Y, z^{(p)}}^2\biggr\}+\mathcal{P} ,
\end{eqnarray*}
where $\mathcal{P}=2PK$ for AIC and $\mathcal{P}=(\log n)PK$ for
BIC, with respect to $K$. Here, for $i\in\mathcal{N}_p$,
$\tilde{\bolds{\epsilon}}_i= \mathbf{V}_i -\tilde{\bolds{\mu}}_{Y,
z^{(p)},i}-\sum_{k=1}^{K}{\tilde{\xi}}_{z^{(p)}, k,
i}^{*}\tilde{\bolds{\phi}}_{z^{(p)}, k, i}$ with $\tilde{\bolds{\mu
}}_{Y, z^{(p)},
i}=({\tilde\mu}_{Y, z^{(p)}}(T_{i1}),
\dots, {\tilde\mu}_{Y, z^{(p)}}(T_{iN_i}))^T$, $\mathbf{V}_i=(V_{i1},
\dots,
V_{iN_i})^T$, $\tilde{\bolds{\phi}}_{z^{(p)}, k, i}=(\tilde\phi_{z^{(p)},
k}(T_{i1}), \dots, \tilde\phi_{z^{(p)}, k}(T_{iN_i}))^T$ and with
estimated principal components
\begin{eqnarray}\label{eq4.1}
{\tilde\xi}_{z^{(p)}, k, i}^*={\tilde\lambda}_{z^{(p)}, k}
\tilde{\bolds{\phi}}_{z^{(p)}, k, i}^T{\tilde\Sigma}_{Y,
z^{(p)},i}^{-1}(\mathbf{V}_i-\tilde{\bolds{\mu}}_{Y,z^{(p)},
i}),
\end{eqnarray}
where $\tilde\Sigma_{Y, z^{(p)},i}$ is an $N_i$-by-$N_i$ matrix whose
$(j,k)$-element is $\tilde G_{Y, z^{(p)}}(T_{ij},
T_{ik})+\tilde\sigma^2_{Y, z^{(p)}}\delta_{jk}$. Analogous criteria
are used for the predictor process $X$, selecting $K$ by minimizing
$\operatorname{AIC}(K)$ and $\operatorname{BIC}(K)$. Marginal
versions of these criteria are also available.

In step [2], for each bin width $h$, we first select the best
smoothing bandwidth $b^{\ast}(h)$ based on AIC or BIC and then
select the final bin width $h^{\ast}$ by a second application of AIC
or BIC, plugging $b^{\ast}(h)$ into this selection as follows. For
a given bin width $h$, define the $P$-by-$P$ smoothing
matrix~$\mathbf{S}_{0,2}$ whose $(p_1, p_2)$th element is $\omega_{0,2}(z^{(p_1)},
z^{(p_2)}, b)$. The effective number of parameters of the
smoothing matrix is then the trace of $\mathbf{S}_{0,2}^T\mathbf{S}_{0,2}$
(cf. \citet{wahba1990}). This suggests minimization of
\begin{eqnarray*}
\operatorname{AIC}(b|h)&=&\sum_{i=1}^{n}\biggl\{\frac{1}{\hat\sigma
^2_{Y}}\hat{\bolds\epsilon}_i^T
\hat{\bolds\epsilon}_i+{N_i}\log(2\uppi)+{N_i}\log
\hat\sigma^2_{Y} \biggr\}+ 2 \operatorname{tr}(\mathbf{S}_{0,2}^T\mathbf{S}
_{0,2}),
\end{eqnarray*}
leading to $b^*(h)$, where
\begin{eqnarray*}
{\hat{\bolds\epsilon}}_i=\mathbf{V}_i-{\hat{\bolds{\mu}}}_{Y,
z_i, i}-\sum_{p=1}^P \omega_{0,2}\bigl(z^{(p)}, z_i, b\bigr) \sum_{m, k=1}^{M,
K}\frac{{\tilde\sigma}_{z^{(p)},mk}}{\tilde\rho_{z^{(p)}, m}}
{\hat\zeta}^*_{z^{(p)},m,i}{\tilde{\bolds{\phi}}}_{z^{(p)},k,i}
\end{eqnarray*}
with ${\hat{\bolds{\mu}}}_{Y, z_i, i}=(\hat
\mu_{Y, z_i}(T_{i1}),
\ldots, \hat\mu_{Y, z_i}(T_{iN_i}))^T$
and estimated principal component scores
\begin{eqnarray*}
{\hat\zeta}_{z^{(p)}, m, i}^*={\tilde\rho}_{z^{(p)}, k}
{\tilde{\bolds{\psi}}}_{z^{(p)}, m, i}^T{\tilde\Sigma}_{X,
z^{(p)},i}^{-1}(\mathbf{U}_i-{\hat{\bolds{\mu}}}_{X,z_i,
i}).
\end{eqnarray*}
The definition of pseudo-BIC scores is analogous.

In step [3], to select the bin width $h^*$, we minimize
\begin{eqnarray*}
\operatorname{AIC}(h,
b^*(h))=\sum_{i=1}^{n}\biggl\{\frac{1}{{\hat\sigma}^2_{Y}}{\hat
{\bolds\epsilon}}_i^T
\hat{\bolds\epsilon}_i+{N_i}\log(2\uppi)+{N_i}\log
\hat\sigma^2_{Y} \biggr\} +2MKP ,
\end{eqnarray*}
or the analogous BIC score, using $b^*(h)$ for each $h$, as determined
in the previous
step.

\section{Simulation study} \label{sec:simulation}

We compare global functional linear regression and
varying-coefficient functional linear regression through simulated
examples with a functional response. For the case of a scalar
response, the proposed varying-coefficient functional linear
regression approach achieves similar performance improvements
(results not reported). For the finite-sample case, there are
several parameters to be selected (see Section \ref{sec:fsi}). In
the simulations, we use pseudo-AIC to select bin width $h$ and
pseudo-BIC to select the smoothing bandwidth $b$ and the number of
regularization terms $M(n)$ and $K(n)$.

The domains of predictor and response trajectories are chosen as
$\mathcal{S}=[0, 10]$ and $\mathcal{T}=[0, 10]$, respectively. The
predictor trajectories $X$ are generated as
$X(s)=\mu_{X}(s)+\sum_{m=1}^{3}\zeta_m \psi_m(s)$ for
$s\in\mathcal{S}$, with mean predictor trajectory
$\mu_{X}(s)=(s+\sin(s))$, the three eigenfunctions are
$\psi_1(s)=-\sqrt{\frac{1}{5}}\cos(\uppi s/5)$,
$\psi_2(s)=\sqrt{\frac{1}{5}}\sin(\uppi s/5)$,
$\psi_3(s)=-\sqrt{\frac{1}{5}}\cos(2\uppi s/5)$ and their
corresponding functional principal components are independently
distributed as $\zeta_1\sim N(0, 2^2)$, $\zeta_2\sim N(0,
\sqrt{2}^2)$, $\zeta_3\sim N(0, 1^2)$. The additional covariate $Z$
is uniformly distributed over $[0, 1]$. For $z\in[0, 1]$, the
slope function is linear in $z$, $\beta(z, s, t)=(z+1)
(\psi_1(s)\psi_1(t)+\psi_2(s)\psi_2(t)+\psi_3(s)\psi_3(t))$ and the
conditional response trajectory is $E(Y(t)|X, Z=z)=\mu_{Y,
z}(t)+\int_{0}^{10} \beta(z, s, t) (X(s)-\mu_{X}(s))\,\mathrm{d}s$, where
$\mu_{Y,z}(t)=(1+z)(t+\sin(t))$. We consider the following two
cases.

\begin{example}[(Regular case)]\label{exfunregular}
The first example focuses on the regular case with dense measurement
design. Observations on the predictor and response trajectories are
made at $s_j=(j-1)/3$ for $j=1, 2, \ldots, 31$ and $t_j=(j-1)/3$ for
$j=1, 2, \ldots, 31$, respectively. We assume the measurement errors
on both the predictor and response trajectories are distributed as
$N(0, 1^2)$, that is, $\sigma_X^2=1$ and $\sigma_Y^2=1$.
\end{example}

\begin{example}[(Sparse and irregular case)]\label{exfunsparse}
In this example, we make a random number of measurements on each
trajectory in the training data set, chosen with equal probability
from $\{2, 3, \ldots, 10\}$. We note that, for the same subject, the
number of measurements on the predictor and the number of
measurements on the response trajectory are independent. For any
trajectory, given the number of measurements, the measurement times
are uniformly distributed over the corresponding trajectory domain.
The measurement error is distributed as $N(0, 1^2)$ for both the
predictor and the response trajectories.
\end{example}

In both examples, the training sample size is 400. An independent
test set of size 1000 is generated with the predictor and response
trajectories fully observed. We compare performance using mean
integrated squared prediction error (MISPE)
\begin{eqnarray*}
&&\frac{1}{1000}\sum_{j=1}^{1000}\int_{\mathcal{T}}\biggl[E(Y_j^*(t)|X_j^*,
Z_j^*)
\\
&&{}\hspace{58pt}-\biggl(\hat\mu_{Y, Z_j^*}(t)+\int_{\mathcal{S}} \hat
\beta(Z_j^*, s, t ) \bigl(X_j^*(s)-\hat\mu_{X, Z_j^*}(s)\bigr)
\,\mathrm{d}s\biggr)\biggr]^2\,\mathrm{d}t\big/|{\mathcal{T}}|,
\end{eqnarray*}
analogously for the global functional linear regression, where
$(X_j^*, Y_j^*, Z_j^*)$ denotes the data of the $j$th subject in
the independent test set. In Table \ref{simulationtable111}, we
report the mean and standard deviation (in parentheses) of the MISPE
of the global and varying-coefficient functional linear
regression over 100 repetitions for each case. This shows that in
this simulation setting, the proposed varying-coefficient functional
linear regression approach reduces MISPE drastically, compared with
the global functional linear regression, both for regular and sparse irregular
designs.

\begin{figure}[b]

\includegraphics{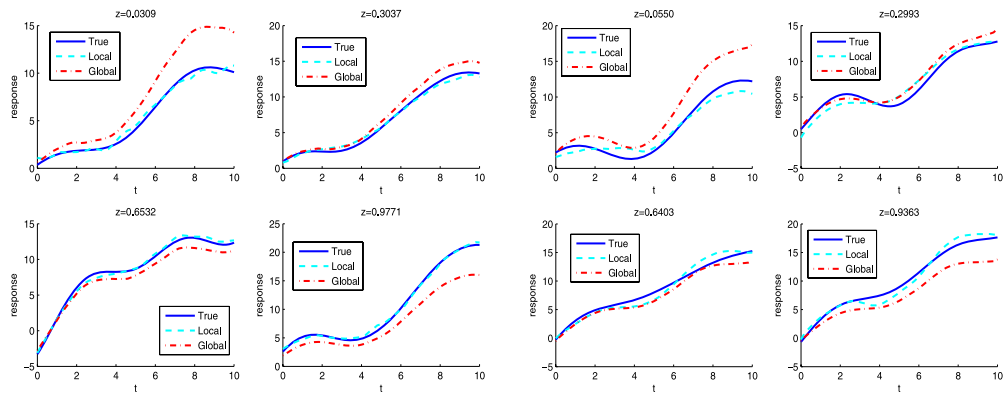}

\caption{In one random repetition, the true (solid) conditional
expected response trajectories and predicted conditional expected
response trajectories via the global functional linear regression
(dot-dashed) and the varying-coefficient functional linear
regression (dashed) are plotted for four randomly selected subjects in
the independent
test set with median integrated squared prediction error. The left four
panels and the right four
correspond to the regular and sparse irregular cases, respectively.}
\label{regularsimupredcomp}
\end{figure}

To visualize the differences between predicted conditional expected
response trajectories, for a small random sample, in both the regular
and sparse and irregular design cases, we randomly choose four
subjects from the test set with median values of the integrated
squared prediction error (ISPE) for the varying-coefficient
functional linear regression. The true and predicted conditional
expected response trajectories are plotted in Figure
\ref{regularsimupredcomp}, where the left four panels correspond to
the regular design case and the right four to the sparse irregular
case. Clearly, the locally varying method is seen to be superior.

\begin{table}[t]
\tablewidth=10.5cm
\caption{Simulation results: mean and standard deviation of MISPE
for global and varying-coefficient functional linear regression with
a functional response, for both regular and sparse
cases}\label{simulationtable111}
\begin{tabular*}{10.5cm}{@{\extracolsep{4in minus 4in}}lll@{}}
\hline
& Functional linear  & Varying-coefficient functional \\
& regression & linear regression\\
\hline
Regular & 4.0146 (1.6115) & 0.7836 (0.4734) \\
Sparse and irregular & 4.0013 (0.8482) & 1.0637 (0.3211)\\
\hline
\end{tabular*}
\end{table}

\section{Applications} \label{sec:application}

We illustrate the comparison of the proposed varying-coefficient
functional linear model with the global functional linear regression in
two applications.

\subsection{Egg-laying data}
The egg-laying data represent the entire reproductive history of one thousand
Mediterranean fruit flies (`medflies' for short), where daily
fecundity, quantified by the number of eggs laid per day, was
recorded for each fly during its lifetime; see \citet{Carey1998}
for details of this data set and experimental background.

We are interested in predicting future egg-laying patterns over an
interval of fixed length, but with potentially different starting
time, based on the daily fecundity information during a fixed
earlier period. The predictor trajectories were chosen as daily
fecundity between day 8 and day~17. This interval covers the tail of
an initial rapid rise to peak egg-laying and the initial part of the
subsequent decline and, generally, the egg-laying behavior at and near
peak egg-laying is included. It is of interest to study in what form
the intensity of peak egg-laying is associated with subsequent
egg-laying behavior, as trade-offs may point to constraints that may
play a role in the evolution of longevity.

\begin{figure}[t]

\includegraphics{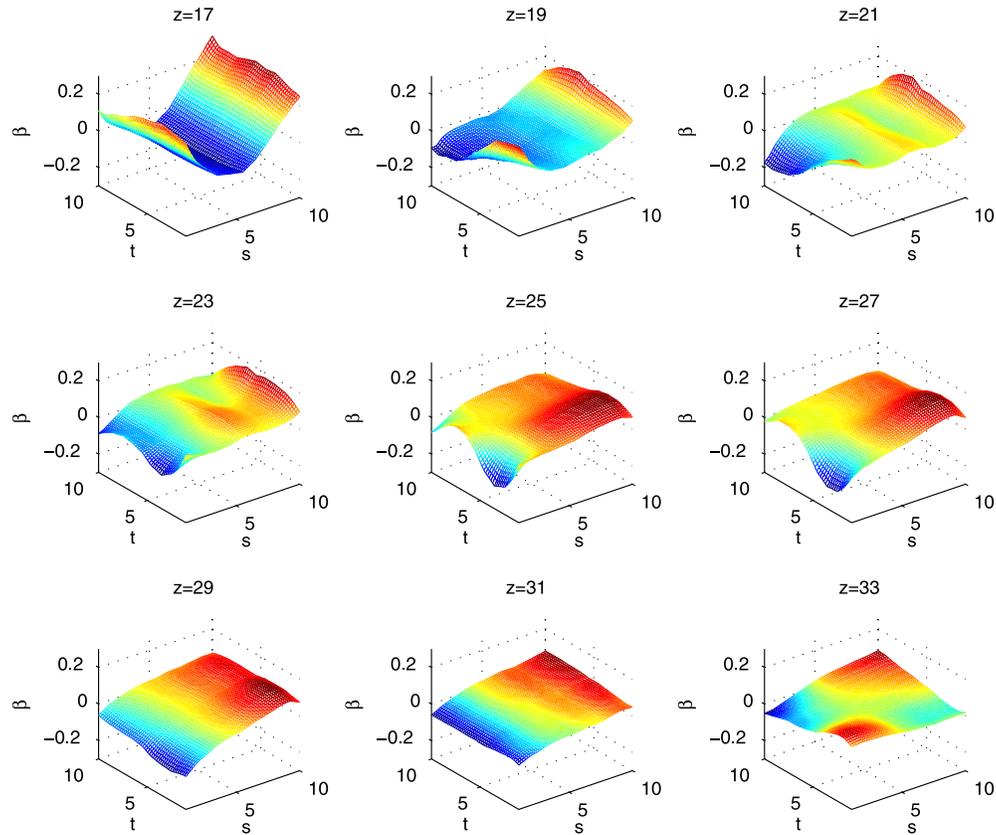}

\caption{The slope functions estimated by the varying-coefficient
functional linear regression at different levels of the additional
covariate $z$ for the egg-laying
data.}\label{egglaylocalbeta}
\end{figure}

\begin{figure}[t]

\includegraphics{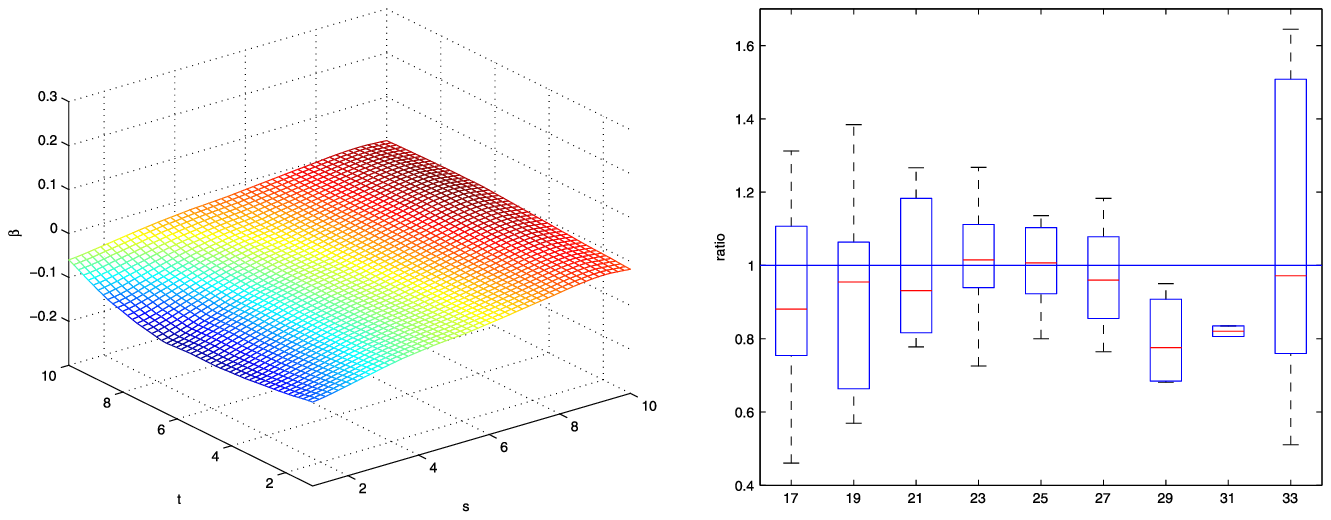}

\caption{The left panel plots the slope function estimated by the
global functional linear regression for the egg-laying data and the
right panel corresponds to box plots of the ratios of MSPE of the
varying-coefficient functional linear regression to that of the
global functional linear regression for the subjects in the test
data set for different levels of the additional covariate~$Z$.}
\label{egglayglobalbeta}
\end{figure}

While the predictor process is chosen with a fixed domain, the
response process has a moving domain, with a fixed length of ten days,
but a different starting age for each subject, which serves as the
additional covariate $Z$. Due to the limited number of subjects in
this study, we use a pre-specified discrete set for the values of
$Z$: $\mathcal{Z}=\{17, 19, 21, 23, 25, 27, 29, 31, 33\}$ with a
pre-specified bin width $h=2$. For subject $i$ with
$z_i\in\mathcal{Z}$, measurements $U_{ij}$ on the predictor
trajectory are the daily numbers of eggs on day $j+7$, and
measurements $V_{ik}$ on the response trajectory correspond to the
daily number of eggs on day $k+z_i$ for $j=1, 2, \dots, 10$ and
$k=1, 2, \dots, 10$. The numbers of subjects in these bins are 30,
29, 18, 29, 22, 19, 19, 17 and 36, respectively. For each bin, we
randomly select 15 subjects as the training set and the remaining
subjects are used to evaluate the prediction performance, comparing
the performance of the global and the varying-coefficient functional
linear regression. The prediction performance is quantified by mean
squared prediction error (MSPE), defined for each subject $i$ in the
test set as
\[
\operatorname{MSPE}_g(i)=\frac{1}{10}\sum_{k=1}^{10}(\hat y_{ik}^g -V_{ik})^2
\quad \mbox{and}\quad  \operatorname{MSPE}_l(i)=\frac{1}{10}\sum_{k=1}^{10}(\hat y_{ik}^l
-V_{ik})^2,
\]
where $\hat y_{ik}^g$ and $\hat y_{ik}^{l}$ denote the predicted
daily fecundities corresponding to $V_{ik}$ using the global (resp.~the
proposed varying-coefficient (local)) functional linear
regression.

Through pseudo-AIC, the global functional linear regression selects
two and three principal components for the predictor and response
trajectories, respectively, while the varying-coefficient functional
linear regression uses two principal components for both trajectories.
After smoothing, the slope functions estimated by the
varying-coefficient models are plotted in Figure
\ref{egglaylocalbeta} for different values of $Z$ and the estimated
slope function for the global functional linear regression is plotted
in the
left panel of Figure \ref{egglayglobalbeta}. Box plots of the ratio
$\operatorname{MSPE}_l(i)/\operatorname{MSPE}_g(i)$ for subjects in the test data
set are shown in the right panel of Figure \ref{egglayglobalbeta}
for different levels of the covariate $Z$. There is one outlier
above the maximum value for $Z=18$ which is not shown. For most
bins, the median ratios are seen to be smaller than 1, indicating
an improvement of our new varying-coefficient functional linear
regression. Denoting the average MSPE (over the independent
test\vadjust{\goodbreak}
data set) of the global and the varying-coefficient functional
linear regression by $\overline{\operatorname{MSPE}}_g$ and
$\overline{\operatorname{MSPE}}_l$, respectively, the relative performance
gain
$(\overline{\operatorname{MSPE}}_l-\overline{\operatorname{MSPE}}_g)/\overline
{\operatorname{MSPE}}_g$
is found to be $-0.0810$ so that the prediction improvement of the
varying-coefficient method is $8.1\%$.

Besides prediction, it is of interest to study the dependency of the
future egg-laying behavior on peak egg-laying. From the changing
slope functions in Figure \ref{egglaylocalbeta}, we find that, for the segments close to
the peak segments, the egg-laying pattern is inverting the peak
pattern, meaning that sharper and higher peaks are associated with
sharp downturns, pointing to a near-future exhaustion effect of peak
egg-laying. In contrast, the shape of egg-laying segments further
into the future is predicted by the behavior of the first derivative
over the predictor segment so that slow declines near the end of
peak egg-laying are harbingers of future robust egg-laying. This is
in accordance with a model of exponential decline in egg-laying
that has been proposed by \citet{MCWLV2001}.

\subsection{BLSA data with scalar response}
As a second example, we use a subset of data from the Baltimore
Longitudinal Study of Aging (BLSA), a major longitudinal data set for
human aging (\citet{Shock1984}, \citet{PMBLF1997}). The data consist of 1590
male volunteers who were scheduled to be seen twice per year.
However, many participants missed scheduled visits or were seen at
other than scheduled times so that the data are sparse and
irregular with unequal numbers of measurements and different
measurement times for each subject. For each subject, current age
and systolic blood pressure (SBP) were recorded during each visit. We
quantify how the SBP trajectories of a subject available in a middle
age range between age 48 and age 53 affect the average of the SBP
measurements made during the last five years included in this study,
at an older age. The predictor domain is therefore of length five years and
the response is scalar. The additional covariate for each subject is
the beginning age of the last five-year interval included in the
study. After excluding subjects with less than two measurements in
the predictor, 214 subjects were included for whom the additional
covariate ranged between 55 and 75. We bin the data according to the
additional covariate, with bin centers at ages 56.0, 59.0, 62.0,
65.0, 68.5 and 73.0 years and the numbers of subjects in each of these
bins are 38, 33, 38, 32, 39 and 34.

\begin{figure}[t]

\includegraphics{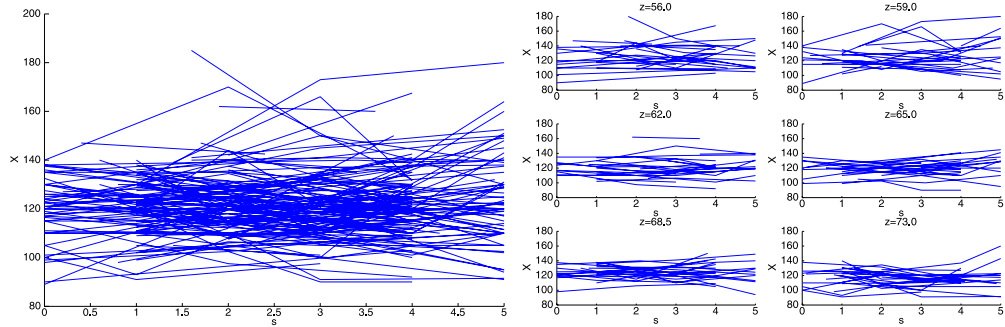}

\caption{Plots of predictor processes: the left panel for the global
functional linear regression and the right panel for different bins
according to the additional covariate in the varying-coefficient
functional linear regression.}\label{scalblsabeta1111}
\end{figure}

\begin{figure}

\includegraphics{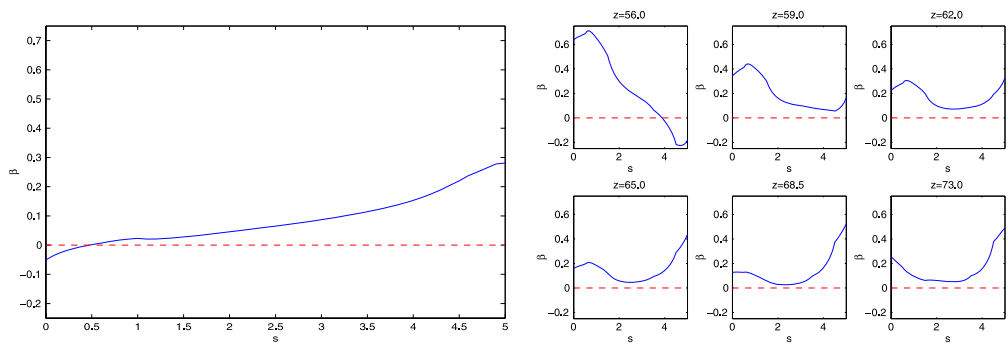}

\caption{The estimated slope function via the global functional
linear regression and the new proposed varying-coefficient
functional linear regression (for different levels of $Z$) are
plotted as the solid lines in the left and right panels,
respectively.}\label{scalblsabeta2222}
\end{figure}

We randomly selected 25 subjects from each bin for model estimation
and used the remaining subjects to evaluate the prediction performance.
In contrast to the egg-laying data, the predictor measurements in
this longitudinal study are sparse and irregular. Pseudo-BIC
selects two principal components for the predictor trajectories for
both global and varying-coefficient functional linear regressions.
Using the same criterion for relative performance gain as in the
previous example, the varying-coefficient functional linear
regression achieves $11.8\%$ improvement compared to the global
functional linear regression. Estimated slope functions are shown in
Figure \ref{scalblsabeta2222} and predictor trajectories in Figure
\ref{scalblsabeta1111}.

The shape changes of the slope functions with changing covariate
indicate that the negative derivative of SBP during the middle-age
period is associated with near-future SBP.  Further into the
future, this pattern is reversed and
an SBP increase
near the right end of the initial period is becoming predictive.

\section{Concluding remarks} \label{sec:conclusion}

Our results indicate that established functional linear regression
models can be improved when an available covariate is incorporated.
We implement this idea by extending the functional linear model to a
varying-coefficient version, inspired by the analogous, highly
successful extension of classical regression models. In both
application examples, the increased flexibility that is inherent in
this extension] leads to clear gains in prediction error. In
addition, it is often of interest to ascertain the effect of the
additional covariate. This can be done by plotting the regression
slopes for each bin defined by the covariate and observing the
dependency of this function or surface on the value of the
covariate.

Further extensions that are of interest in many applications concern
the case of multivariate covariates. If the dimension is low, the
smoothing methods and binning methods that we propose here can be
extended to this case. For higher-dimensional covariates or
covariates that are not continuous, one could form a single index to
summarize the covariates and thus create a new one-dimensional
covariate which then enters the functional regression model in the
same way as the one-dimensional covariate that we consider.

As seen in the data applications, the major applications of the
proposed methodology are expected to come from longitudinal studies
with sparse and irregular measurements, where the presence of
additional non-functional covariates is common.

\appendix
\section*{Appendix: Auxiliary results and proofs}
\label{app}

We note that further details, such as omitted proofs, can be
found in a technical report that is available at
\url{http://www4.stat.ncsu.edu/\textasciitilde wu/WuFanMueller.pdf}.

A bivariate kernel function
$\kappa_2(\cdot, \cdot)$ is said to be of order $(\nu, \ell)$ with
$\nu=(\nu_1, \nu_2)$ if it satisfies
%
\begin{eqnarray}\label{2dkernelorder}
\int u^{\ell_1}v^{\ell_2} \kappa_2(u,v)\, \mathrm{d}u\,\mathrm{d}v
=\cases{ 0, &\quad  $0 \leq
\ell_1+\ell_2 < \ell, \ell_1\ne\nu_1, \ell_2\ne\nu_2,$
\cr
\nu!, &\quad $\ell_1=\nu_1, \ell_2=\nu_2,$
\cr
\ne0, &\quad $\ell_1+\ell_2=\ell$,
}
\end{eqnarray}
and
%
\begin{eqnarray}\label{2dkernelorder2}
\int|u^{\ell_1}v^{\ell_2} \kappa_2(u,v)| \,\mathrm{d}u\,\mathrm{d}v <\infty
\qquad \mbox{for any } \ell_1+\ell_2=\ell,
\end{eqnarray}
where $\nu!=\nu_1!\cdot\nu_2!$. Similarly, a univariate kernel
function $\kappa_1(\cdot)$ is of order $(\nu, \ell)$ for a
univariate $\nu=\nu_1$ when (\ref{2dkernelorder}) and
(\ref{2dkernelorder2}) hold for $\ell_2\equiv0$ on the right-hand
side while integrating over the univariate argument $u$ on the left.

We introduce the following technical conditions:
\begin{enumerate}[(viii)]
\item[(i)] The variable $S$ has compact domain
$\mathcal{S}$. Given $Z=z$, $S$ has conditional density
$f_{S, z}(s)$. Assume, uniformly in
$z\in\mathcal{Z}$, that $\frac{\partial^\ell}{\partial
s^\ell}f_{S,z}(s)$ exists and is continuous for $\ell=2$\vspace*{1pt} on
$\mathcal{S}$ and, further,
$\inf_{s\in\mathcal{S}} f_{S,z}(s)>0$, analogously for $T$.

\item[(ii)] Denote the conditional density
functions of $(S, U)$ and $(T, V)$
by $g_{X, z}(s, u)$ and $g_{Y, z}(t, v)$,
respectively. Assume that the derivative $\tfrac{\partial
^\ell}{\partial
s^\ell}g_{X,z}(s, u)$
exists for all\vspace{1pt} arguments $(s,
u)$, is uniformly continuous on $\mathcal{S}\times
\mathbb{R}$ and is Lipschitz continuous
in $z$, for $\ell=2$, analogously for $g_{Y,z}(t, v)$.

\item[(iii)] Denote the conditional density
functions of quadruples $(S_1, S_2, U_1, U_2)$ and $(T_1, T_2,
V_1, V_2)$ by $g_{2X, z}(s_1, s_2, u_1, u_2)$ and $g_{2Y, z}(t_1,
t_2, v_1, v_2)$, respectively. For simplicity, the corresponding marginal
conditional densities of $(S_1, S_2)$ and $(T_1, T_2)$ are also
denoted by $g_{2X,z}(s_1, s_2)$ and $g_{2Y,z}(t_1, t_2)$,
respectively. Denote the conditional density of $(S, T, U, V)$ given
$Z=z$ by $g_{XY, z}(s, t, u, v)$ and, similarly, its corresponding
conditional marginal density of $(S, T)$ by $g_{XY,z}(s, t)$.
Assume that the derivatives
$\frac{\partial^\ell}{\partial s_1^{\ell_1}\,\partial
s_2^{\ell_2}} g_{2X,z}(s_1, s_2, u_1, u_2)$
exist for all arguments
$(s_1, s_2, u_1, u_2)$, are uniformly
continuous on $\mathcal{S}^2\times\mathbb{R}^2$ and are Lipschitz
continuous in $z$ for
$\ell_1+\ell_2=\ell$, $0\leq\ell_1, \ell_2\leq\ell=2$,
analogously for $g_{2Y,z}(t_1, t_2, v_1, v_2)$ and $g_{XY,z}(s, t, u, v)$.

\item[(iv)] For every $p=1, 2, \dots, P$,
$b_{X,z^{(p)}}\rightarrow0$, $n_{z^{(p)},h}
b_{X,z^{(p)}}^{4}\rightarrow\infty$, $n_{z^{(p)},h}
b_{X,z^{(p)}}^{6}<\infty$, $b_{Y,z^{(p)}}\rightarrow0$,
$n_{z^{(p)},h} b_{Y,z^{(p)}}^{4}\rightarrow\infty$ and
$n_{z^{(p)},h} b_{Y,z^{(p)}}^{6}<\infty$ as $n\rightarrow\infty$.

\item[(v)] For every $p=1, 2, \dots, P$,
$h_{X,z^{(p)}}\rightarrow0$, $n_{z^{(p)},h}
h_{X,z^{(p)}}^{6}\rightarrow\infty$, $n_{z^{(p)},h}
h_{X,z^{(p)}}^{8}<\infty$, $h_{Y,z^{(p)}}\rightarrow0$,
$n_{z^{(p)},h} h_{Y,z^{(p)}}^{6}\rightarrow\infty$ and
$n_{z^{(p)},h} h_{Y,z^{(p)}}^{8}<\infty$ as $n\rightarrow\infty$.

\item[(vi)] For every $p=1, 2, \dots, P$,
$h_{1,z^{(p)}}/h_{2,z^{(p)}}\rightarrow1$, $h_{1,z^{(p)}}\rightarrow0$,
$n_{z^{(p)},h} h_{1,z^{(p)}}^{6}\rightarrow\infty$ and $n_{z^{(p)},h}
h_{1,z^{(p)}}^{8}<\infty$ as $n\rightarrow\infty$.\vspace{1pt}

\item[(vii)] For every $p=1, 2, \dots, P$,
$b_{X,z^{(p)}, V}\rightarrow0$, $n_{z^{(p)},h}
b_{X,z^{(p)}, V}^{4}\rightarrow\infty$, $n_{z^{(p)},h} b_{X,z^{(p)},
V}^{6}<\infty$, $b_{Y,z^{(p)}, V}\rightarrow0$, $n_{z^{(p)},h}
b_{Y,z^{(p)}, V}^{4}\rightarrow\infty$ and $n_{z^{(p)},h}
b_{Y,z^{(p)}, V}^{6}<\infty$ as $n\rightarrow\infty$.\vspace{1pt}

\item[(viii)] Univariate kernel $\kappa_1$ and
bivariate kernel
$\kappa_2$ are compactly supported, absolutely integrable and of
orders $(\nu,\ell)=(0, 2)$ and $((0, 0), 2)$, respectively.

\item[(ix)] Assume that $\sup_{(z,
s)\in\mathcal{Z}\times\mathcal{S}}
E(E(X(s)-\mu_{X,Z}(s))^4|Z=z))<\infty$, and analogously for~$Y$.

\item[(x)] The slope function $\beta(z, s, t)$ is
twice differentiable in $z$, that is, for any $(s, t)\in\mathcal
{S}\times\mathcal{T}$, $\frac{\partial^2}{\partial z^2}\beta(z, s, t)$
exists and is continuous in $z$.\vspace{1,5pt}

\item[(xi)] The bin width $h$ and smoothing bandwidth $b$
are such that $b/h<\infty$ as $n\rightarrow\infty$. The bin width
$h$ is chosen such that
$P\propto n^{1/8}$.
\end{enumerate}

\begin{prop}\label{aprobprop}
For $E_n$ defined in (\ref{aprobdefn}), under \textup{(xi)}, it
holds that $P(E_n)\rightarrow1$ as $n\rightarrow\infty$.
\end{prop}

\begin{pf} First, note that
$P(\min n_{z^{(p)}, h}>\tilde n)\geq1-\sum_{p=1}^{P} P(n_{z^{(p)},
h}< \tilde n).$
Consider the $p$th bin and let $\pi_p=P(Z\in[z^{(p)}-\frac{h}{2},
z^{(p)}-\frac{h}{2}))$. Then $n_{z^{(p)}, h}$ is asymptotically
distributed as $N(n\pi_p, n\pi_p(1-\pi_p))$ due to the normal
approximation to a binomial random variable. Thus, $P( n_{z^{(p)}, h} >
\tilde n)\rightarrow f_{N(0, 1)}(a_p)/a_p$ with $a_p=-(\tilde n - n\pi
_p)/\sqrt{n \pi_p(1-\pi_p)}$, where $f_{N(0, 1)}(\cdot)$ is the
probability density function of the standard normal distribution.
Due to \textup{[A1]}, $\pi_p$ is bounded between $\underline
{f_{Z}}/(\underline{f_{Z}} + (P-1)\bar{f_Z})$ and $\bar
{f_Z}/((P-1)\underline{f_{Z}} + \bar{f_Z})$. It follows that
$P(E_n)\rightarrow1$ as $n\rightarrow\infty$ by noting that $\tilde
n \propto\sqrt{n}$, $P\propto n^{1/8}$, and $f_{N(0, 1)}(x)/x$ decays
exponentially in $x$.
\end{pf}

We next prove the consistency of the raw estimate of the mean functions
of predictor and response trajectories within each bin. Consider a
generic bin $[z-h/2, z+h/2)$, with bin center $z$ and bandwidth
$h$, and
let $b_{X,z}$ and $b_{Y,z}$ be smoothing bandwidths used to
estimate $\mu_{X,z}(s)$ and $\mu_{Y,z}(t)$, $h_{X,z}$ and $h_{Y,z}$
for $G_{X,z}(s_1, s_2)$ and $G_{Y,z}(t_1, t_2)$, respectively,
$h_{1,z}$ and
$h_{2,z}$ for $C_{XY,z}(s,t)$, and $b_{X, z, V}$ and $b_{Y,z, V}$ for
$V_{X,z}(s)=G_{X,z}(s, s)+\sigma_X^2$ and
$V_{Y,z}(t)=G_{Y,z}(t,t)+\sigma_Y^2$, respectively.

For a positive integer $l\geq1$, let $\{\psi_{p}(t, v), p=1, 2,
\dots, l\}$ be a collection of real functions $\psi_p\dvtx
\mathbb{R}^2\rightarrow\mathbb{R}$ satisfying the following
conditions:
\begin{enumerate}[{[C1.1a]}]
\item[{[C1.1a]}] The derivative functions $\frac
{\partial^\ell}{\partial t^\ell}\psi_p(t,
v)$ exist for all arguments $(t, v)$ and are uniformly continuous on
$\mathcal{T}\times\mathbb{R}$.
\item[{[C1.2a]}] $\int\int\psi_p^2(t, v) g_{Y, z}(t, v)\,\mathrm{d}v \, \mathrm{d}t <
\infty$.
\item[{[C2.1a]}] Uniformly in $z\in\mathcal{Z}$, bandwidths
$b_{Y,z}$ for one-dimensional smoothers
are such that $b_{Y,z}\rightarrow0$, $n_{z,h}
b_{Y,z}^{\nu+1}\rightarrow\infty$ and $n_{z,h}
b_{Y,z}^{2\ell+2}<\infty$ as $n\rightarrow\infty$.
\end{enumerate}

Define $\mu_{p\psi, z}=\mu_{p\psi, z}(t)=\frac{\mathrm{d}^\nu}{\mathrm{d}t^\nu
}\int\psi_p(t,
v)g_{Y, z}(t, v)\,\mathrm{d}v$ and
\begin{eqnarray*}
\Psi_{pn,z}=\Psi_{pn,z}(t)=\frac{1}{n_{z, h}
b_{Y,z}^{\nu+1}}\sum_{i\in
\mathcal{N}_{z,h}}\frac{1}{EN}\sum_{j=1}^{N_i}\psi_p(T_{ij},
V_{ij})\kappa_1\biggl(\frac{T_{ij}-t}{b_{Y,z}}\biggr),
\end{eqnarray*}
where $g_{Y, z}(t, v)$ is the conditional density of $(T, V)$, given
$Z=z$.

\begin{lem} \label{ymwjasalemma1} Under conditions \textup{[A0]--[A3]}
 \textup{(i)},
\textup{(ii)}, \textup{(viii)},
\textup{[C1.1a]}, \textup{[C1.2a]} and \textup{[C2.1a]}, we have $\tau_{pn}=\sup_{(z,t)\in
\mathcal{Z}\times\mathcal{T}} {\vert
\Psi_{pn,z}(t)-\mu_{p\psi, z}(t)\vert}/({h+ (\sqrt{n_{z,h}}
b_{Y,z}^{\nu+1})^{-1}})=\mathrm{O}_p(1).$
\end{lem}

\begin{pf} Note that
$\vert\Psi_{pn,z}(t)-\mu_{p\psi, z}(t)\vert\leq\vert\Psi_{pn,z}(t)-
E \Psi_{pn,z}(t)\vert+ \vert E\Psi_{pn,z}(t)- \mu_{p\psi,
z}(t)\vert$
and $E\vert\tau_{pn}
\vert=\mathrm{O}(1)$ implies that $\tau_{pn}= \mathrm{O}_p(1)$. Standard conditioning
techniques lead to
\begin{eqnarray*}
E\Psi_{pn,z}(t) &=&\frac{1}{b_{Y,z}^{\nu+1}}E\biggl(E\biggl(\psi_p(T_{i1},
V_{i1})\kappa_1\biggl(\frac{T_{i1}-t}{b_{Y,z}}\biggr)\Bigm| z-\frac{h}{2}\leq Z_i
< \frac{h}{2}\biggr)\biggr).
\end{eqnarray*}
For $Z_i=z_i\in[z-h/2, z+h/2)$, perform a Taylor expansion of
order $\ell$ on the integrand:
\begin{eqnarray*}
&&E\biggl[\psi_p(T_{i1},
V_{i1})\kappa_1\biggl(\frac{T_{i1}-t}{b_{Y,z}}\biggr)\biggr]
\\
&&\quad = \int\int\psi_p(t_1,
v_1) g_{Y,z_i}(t_1, v_1)
\kappa_1\biggl(\frac{t_1-t}{b_{Y,z}}\biggr)\, \mathrm{d} t_1 \,\mathrm{d}v_1
\\
&&\quad =\int\int\biggl( \frac{\partial^\nu}{\partial t^\nu} (\psi_p(t,
v_1) g_{Y,z_i}(t, v_1))\biggr)
\frac{(t_1-t)^\nu}{\nu!}\kappa_1\biggl(\frac{t_1-t}{b_{Y,z}}\biggr) \,\mathrm{d} t_1
\,\mathrm{d}v_1
\\
&& {}\qquad + \int
\int\biggl(\frac{\partial^\ell}{\partial t^\ell} (\psi_p(t, v_1)
g_{Y,z_i}(t, v_1))\biggr)\biggm|_{t=t^*} \frac{ (t_1-t)^{\ell}}{\ell!}
\kappa_1\biggl(\frac{t_1-t}{b_{Y,z}}\biggr)\, \mathrm{d} t_1\, \mathrm{d}v_1,
\end{eqnarray*}
where $t^*$ is between $t$ and $t_1$. Hence,
$|E[\psi_p(T_{i1},
V_{i1})\kappa_1(\frac{T_{i1}-t}{b_{Y,z}})]-\mu_{p\psi, z_i}(t)
b_{Y,z}^{\nu+1}|\leq c_0 \frac{b_{Y,z}^{\ell+1}}{\ell!}
\times\break \int| u^\ell\kappa_1(u)|\,\mathrm{d}u$
due to [C1.2a] and the assumption that the kernel function
$\kappa_1(\cdot)$ is of type $(\nu, \ell)$, where $c_0$ is bounded
according to [C1.1a],
$c_0\leq\sup_{(z_i, t)\in
\mathcal{Z}\times\mathcal{T}}|\frac{\partial^\ell}{\partial
t^\ell} \int\psi_p(t, v_1)\* g_{Y,z_i}(t, v_1)\,\mathrm{d}v_1| < \infty.$
Furthermore, using (ii), we may bound
\begin{eqnarray}\label{lemmabound1}
&&\sup_{t\in\mathcal{T}}\vert E\Psi_{pn,z}(t)-\mu_{p\psi,z}(t)\vert\nonumber
\\
&&\quad \leq c_0 {b_{Y,z}^{\ell-\nu}}/{(\ell!)} \int| u^\ell
\kappa_1(u)|\,\mathrm{d}u
\\
&&{}\qquad +E\biggl\{E\biggl[\sup_{t\in\mathcal{T}}
|\mu_{p\psi,
Z_i}(t)-
\mu_{p\psi, z}(t)|\Bigm| z-\frac{h}{2}\leq Z_i < \frac
{h}{2}\biggr]\biggr\}\nonumber
\\
&&\quad \leq c_0 \biggl(\int| u^\ell\kappa_1(u)|\,\mathrm{d}u\biggr)
b_{Y,z}^{\ell-\nu}/(\ell!)+c_1h, \nonumber
\end{eqnarray}
where the constants do not depend on $z$. To bound $E\sup_{t\in
\mathcal{T}}\vert\Psi_{pn,z}(t)- E \Psi_{pn,z}(t)\vert$, we denote
the Fourier transform of $\kappa_1(\cdot)$ by $\zeta_1(t)=\int
\mathrm{e}^{-\mathrm{i}ut}\kappa_1(u)\,\mathrm{d}u$, and letting $\varphi_{pn,
z}(u)=\frac{1}{n_{z,h}} \sum_{m\in\mathcal{N}_{z, h}} \frac{1}{EN}
\sum_{j=1}^{N_m} \mathrm{e}^{\mathrm{i}uT_{mj}}\psi_p(T_{mj}, Y_{mj})$, we have
\begin{eqnarray*}
\Psi_{pn, z}&=&\frac{1}{n_{z,h} b_{Y,z}^{\nu+1}} \sum_{m\in
\mathcal{N}_{z,h}}\frac{1}{EN} \sum_{j=1}^{N_m}
\kappa_1\biggl(\frac{T_{mj}-t}{b_{Y,z}}\biggr)\psi_p(T_{mj},
Y_{mj})\nonumber
\\
&=& \frac{1}{2\uppi b_{Y,z}^{\nu}}\int\varphi_{pn,
z}(u)\mathrm{e}^{-\mathrm{i}tu}\zeta_1(u b_{Y,z}) \,\mathrm{d}u\nonumber
\end{eqnarray*}
and
$\sup_{t\in\mathcal{T}}\vert\Psi_{pn,z}(t)- E \Psi_{pn,z}(t)\vert
\leq\frac{1}{2\uppi b_{Y,z}^\nu}\int\vert\varphi_{pn, z}(u)- E
\varphi_{pn, z}(u)\vert\cdot\vert\zeta_1(u b_{Y,z})\vert
\,\mathrm{d}u.$

Decomposing $\varphi_{pn,z}(\cdot)$ into real and imaginary parts,
\begin{eqnarray*}
\varphi_{pn, z, R}(u)&=&\frac{1}{n_{z,h}} \sum_{m\in\mathcal{N}_{z,
h}} \frac{1}{EN} \sum_{j=1}^{N_m} \cos{(uT_{mj})}\psi_p(T_{mj},
Y_{mj}) ,
\\
\varphi_{pn, z, I}(u)&=&\frac{1}{n_{z,h}} \sum_{m\in\mathcal{N}_{z,
h}} \frac{1}{EN} \sum_{j=1}^{N_m} \sin{(uT_{mj})}\psi_p(T_{mj},
Y_{mj}),
\end{eqnarray*}
we obtain $E \vert\varphi_{pn, z}(u)- E \varphi_{pn, z}(u)\vert=E
\vert
\varphi_{pn, z,R}(u)- E \varphi_{pn, z, R}(u)\vert+\break
E \vert \varphi_{pn, z,I}(u) - E \varphi_{pn, z, I}\! (u)\vert$. Note the
inequality $E \vert\varphi_{pn, z, R}(u)- E \varphi_{pn, z, R}(u)\vert
\leq\break \sqrt{\!E \vert\varphi_{pn, z, R}(u)- E \varphi_{pn, z,
R}(u)\vert^2\!}$ and the fact that $\{[Z_i, N_i, (T_{ij},
Y_{ij})_{j=1}^{N_i}]\dvt i\in
\mathcal{N}_{z, h}\}$ are i.i.d.~implies that
\begin{eqnarray*}
\operatorname{var}(\varphi_{pn, z, R}(u))
\leq\frac{1}{n_{z,h}}E\bigl\{ E\bigl(\psi_p^2(T_{m1}, Y_{m1})|z-h/2\leq Z_m
< z+h/2\bigr)\bigr\} ,
\end{eqnarray*}
where $m\in N_{z,h}$, analogously for the imaginary part.
As a result, we have
\begin{eqnarray*}
&&E\sup_{t\in\mathcal{T}}\vert\Psi_{pn,z}(t)- E \Psi_{pn,z}(t)\vert
\\
&&\quad \leq\frac{2\sqrt{E\{ E(\psi_p^2(T_{m1}, Y_{m1})|z-h/2\leq Z_m <
z+h/2)\}}\int\vert\zeta_1(u)\vert\, \mathrm{d}u }{2\uppi\sqrt{n_{z,h}}
b_{Y,z}^{\nu+1}}.
\end{eqnarray*}
Note that $E(\psi_p^2(T_{m1}, Y_{m1}))$ as a function of $Z_m$ is
continuous over the compact domain $\mathcal{Z}$ and is consequently bounded.
Let
$c_2=2\sup_{Z_m\in\mathcal{Z}}\sqrt{E(\psi_p^2(T_{m1},
Y_{m1}))}<\infty$. Hence, we have
%
\begin{equation}\label{lemmabound2}
E \sup_{t\in\mathcal{T}} \vert\Psi_{pn,z}(t)- E \Psi_{pn,z}(t)\vert
\leq\frac{c_2\int|\zeta_1(u)|\,\mathrm{d}u}{2\uppi} \bigl(\sqrt{n_{z,h}}
b_{Y,z}^{\nu+1}\bigr)^{-1},
\end{equation}
where the constant $c_2(\int|\zeta_1(u)|\,\mathrm{d}u)/{(2\uppi)}$
does not depend on $z$.

The result follows as condition [A1] implies that
$n_{z,h}$ goes to infinity uniformly for $z\in\mathcal{Z}$ as
$n\rightarrow\infty$ and $n_{z,h} b_{Y,z}^{2\ell+2}< \infty$ implies
that $b_{Y,z}^{\ell-\nu}=\mathrm{O}(1/(\sqrt{n_{z,h}} b_{Y,z}^{\nu+1}))$.
We next extend\vspace{1pt} Theorem 1
in \citet{YMW2005jasa} under some additional conditions.
\end{pf}

\begin{enumerate}[{[C3]}]
\item[{[C3]}] Uniformly in $z\in\mathcal{Z}$, $b_{X,z}\rightarrow
0$, $n_{z,h}
b_{X,z}^{4}\rightarrow\infty$, $n_{z,h} b_{X,z}^{6}<\infty$,
$b_{Y,z}\rightarrow0$, $n_{z,h} b_{Y,z}^{4}\rightarrow\infty$ and
$n_{z,h} b_{Y,z}^{6}<\infty$ as $n\rightarrow\infty$.
\end{enumerate}

\begin{lem}\label{ymwjasathm11dmean}
Under conditions \textup{[A0]--[A3]}, (\textup{i), (ii),
(viii), (ix)}
and \textup{[C3]}, we have
\begin{eqnarray}\label{meanyrawcons}
\sup_{(z,s)\in\mathcal{Z}\times\mathcal{S}}\frac{\vert\tilde\mu_{X,
z}(s)- \mu_{X, z}(s)\vert}{h+(\sqrt{n_{z,h}}
b_{X,z})^{-1}}&=&\mathrm{O}_p(1) \quad \mbox{and}\nonumber
\\[-8pt]\\[-8pt]
\sup_{(z,t)\in\mathcal{Z}\times\mathcal{T}}\frac{\vert\tilde\mu_{Y,
z}(t)- \mu_{Y, z}(t)\vert}{h+(\sqrt{n_{z,h}}
b_{Y,z})^{-1}}&=&\mathrm{O}_p(1).\nonumber
\end{eqnarray}
\end{lem}

\begin{pf}
The proof is similar to the proof of Theorem 1
in \citet{YMW2005jasa}.
\end{pf}

Our next two lemmas concern the consistency for estimating the covariance
functions, based on the
observations in the generic bin $[z-h/2, z+h/2)$.
Let $\{\theta_{p}(r_1, r_2, v_1, v_2),
p=1, 2, \dots, l\}$ be a collection of real functions $\theta_p\dvtx
\mathbb{R}^4\rightarrow\mathbb{R}$ with the following properties:

\begin{enumerate}[{[C1.1b]}]
\item[{[C1.1b]}] the derivatives
$\frac{\partial^\ell}{\partial r_1^{\ell_1}\,\partial
r_2^{\ell_2}}\theta_p(r_1, r_2, v_1, v_2)$ exist for all arguments
$(r_1, r_2, v_1, v_2)$ and\vspace{-2pt} are uniformly continuous on
$\mathcal{R}_1\times\mathcal{R}_2\times\mathbb{R}^2$ for
$\ell_1+\ell_2=\ell$, $0\leq\ell_1, \ell_2\leq\ell$, $\ell=0,
1,2$;

\item[{[C1.2b]}] the expectation $\int\int\int\int
\theta_p^2(r_1,
r_2, v_1, v_2) g(r_1, r_2, v_1, v_2)\,\mathrm{d}r_1\,\mathrm{d}r_2\,\mathrm{d}v_1\,\mathrm{d}v_2$ exists and is
finite, uniformly bounded on $\mathcal{Z}$;

\item[{[C2.1b]}] uniformly in $z\in\mathcal{Z}$, bandwidths
$h_{Y,z}$ for the two-dimensional
smoother satisfy $h_{Y,z}\rightarrow0$,
$n_{z,h}h_{Y,z}^{|\nu|+2}\rightarrow\infty$,
$n_{z,h}h_{Y,z}^{2\ell+4}<\infty$ as $n\rightarrow\infty$.
\end{enumerate}

Define $\varrho_{p\theta, z}=\varrho_{p\theta, z}(t_1,
t_2)=\frac{\partial^{|\nu|}}{\partial t_1^{\nu_1}\, \partial
t_2^{\nu_2}} \int\int\theta_p(t_1, t_2, v_1, v_2)g_{2Y,z}(t_1, t_2,
v_1, v_2)\,\mathrm{d}v_1 \,\mathrm{d}v_2$ and
\begin{eqnarray*}
&&\Theta_{pn,z}(t_1, t_2)
=\frac{1}{n_{z,h} h_{Y,z}^{|\nu|+2}}\sum_{i\in\mathcal{N}_{z,h}}\frac{1}{EN(EN-1)}
\\
&&{}\hspace{127pt}\times \sum_{1\leq j\ne
k\leq N_i}\theta_p(T_{ij}, T_{ik}, V_{ij}, V_{ik})
\kappa_2\biggl(\frac{T_{ij}-t_1}{h_{Y,z}},
\frac{T_{ik}-t_2}{h_{Y,z}}\biggr).
\end{eqnarray*}

\begin{lem} \label{ymwjasalemma12dvar}
Under conditions \textup{[A0]--[A3]}, \textup{(i),
(ii), (iii),
(viii), [C1.1b]} with $\mathcal{R}_1=\mathcal{T}$ and
$\mathcal{R}_2=\mathcal{T}$, \textup{[C1.2b]} with $g(\cdot, \cdot, \cdot,
\cdot)=g_{2Y,z}(\cdot, \cdot, \cdot, \cdot)$ and \textup{[C2.1b]}, we have
\[
\vartheta_{pn}=\sup_{(z, t_1, t_2)\in
\mathcal{Z}\times\mathcal{T}\times\mathcal{T}}
\frac{\vert\Theta_{pn,z}-\varrho_{p\theta,z}\vert}{h+(\sqrt{n_{z,h}}
h_{Y,z}^{|\nu|+2})^{-1}}=\mathrm{O}_p(1).
\]
\end{lem}

\begin{pf}
This is analogous to the proof of Lemma 1.
\end{pf}

\begin{enumerate}[{[C4]}]
\item[{[C4]}] Uniformly in $z\in\mathcal{Z}$, $h_{X,z}\rightarrow
0$, $n_{z,h}
h_{X,z}^{6}\rightarrow\infty$, $n_{z,h} h_{X,z}^{8}<\infty$,
$h_{Y,z}\rightarrow0$, $n_{z,h} h_{Y,z}^{6}\rightarrow\infty$ and
$n_{z,h} h_{Y,z}^{8}<\infty$ as $n\rightarrow\infty$.
\end{enumerate}

The proof of the next result is omitted.
\begin{lem} \label{ymwjasathm12dvar} 
Under conditions \textup{[A0]--[A3], (i)--(iii),
(viii), (ix), [C3]} and \textup{[C4]}, we
have
%
\begin{eqnarray}\label{covxunifconrate}
\sup_{(z, s_1, s_2)\in
\mathcal{Z}\times\mathcal{S}^2}\frac{\vert\tilde G_{X, z}(s_1, s_2)-
G_{X, z}(s_1, s_2)\vert}{(h+(\sqrt{n_{z,h}}
h_{X,z}^2)^{-1})}&=&\mathrm{O}_p(1) ,
\\\label{covyunifconrate}
\sup_{(z,t_1, t_2)\in
\mathcal{Z}\times\mathcal{T}^2}\frac{\vert\tilde G_{Y, z}(t_1, t_2)-
G_{Y, z}(t_1, t_2)\vert}{(h+(\sqrt{n_{z,h}}
h_{Y,z}^2)^{-1})}&=&\mathrm{O}_p(1) .
\end{eqnarray}
\end{lem}

To estimate variance of the measurement errors, as in
\citet{YMW2005jasa}, we first estimate $G_{X,z}(s, s)+\sigma^2_X$
(resp.~$G_{Y,z}(t,t)+\sigma^2_Y$) using a local linear smoother
based on $G_{X, i, z}(S_{il}, S_{il})$ for $l=1, 2, \dots, L_i$,
$i\in\mathcal{N}_{z,h}$ (resp.~$G_{Y, i, z}(T_{ij}, T_{ij})$ for
$j=1, 2, \dots, N_i$, $i\in\mathcal{N}_{z,h}$) with smoothing
bandwidth $b_{X,z,V}$ (resp.~$b_{Y,z,V}$) and denote the
estimates\vspace{1pt}
by $\tilde V_{X,z}(s)$ (resp.~$\tilde V_{Y,z}(t)$), removing the two
ends of the interval $\mathcal{S}$ (resp.~$\mathcal{T}$) to get more
stable estimates of $\sigma^2_X$ (resp.~$\sigma^2_Y$). Denote the
estimates based on the generic bin $[z-h/2, z+h/2)$ by $\tilde
\sigma_{X,z}^2$ and $\tilde\sigma_{Y,z}^2$, let
$|\mathcal{S}|$ denote the length of $\mathcal{S}$ and let
$\mathcal{S}_1=[\inf\{s\dvt
s\in\mathcal{S}\}+|\mathcal{S}|/4, \sup\{s\dvtx
s\in\mathcal{S}\}-|\mathcal{S}|/4]$. Then
\begin{eqnarray}
\tilde\sigma_{X,z}^2=\frac{2}{|\mathcal{S}|}\int
_{\mathcal{S}_1}
[\tilde V_{X}(s) -\tilde G_{X,z}(s,s)]\,\mathrm{d}s ,\nonumber
\end{eqnarray}
and analogously for $\tilde\sigma_{Y,z}^2$. Lemmas
\ref{ymwjasathm11dmean} and \ref{ymwjasathm12dvar} imply the
convergence of $\tilde\sigma_{X,z}^2$ and $\tilde\sigma_{Y,z}^2$,
as stated in Corollary \ref{sigmacon}.

\begin{enumerate}[{[C5]}]
\item[{[C5]}] Uniformly in $z\in\mathcal{Z}$, $b_{X,z,V}\rightarrow
0$, $n_{z,h}
b_{X,z,V}^{4}\rightarrow\infty$, $n_{z,h} b_{X,z,V}^{6}<\infty$,
$b_{Y,z,V}\rightarrow0$, $n_{z,h} b_{Y,z,V}^{4}\rightarrow\infty$
and $n_{z,h} b_{Y,z,V}^{6}<\infty$ as $n\rightarrow\infty$.
\end{enumerate}

\begin{coro}\label{sigmacon} Under condition \textup{[C5]} and the conditions
of Lemmas \ref{ymwjasathm11dmean}
and \ref{ymwjasathm12dvar},
\begin{eqnarray*}
\sup_{z\in\mathcal{Z}} {|\tilde\sigma_{X,z}^2-
\sigma_X^2|}/{ \bigl(h+ \bigl(\sqrt{n_{z,h}} b_{X,z, V}\bigr)^{-1}
+\bigl(\sqrt{n_{z,h}}
h_{X,z}^2\bigr)^{-1}\bigr)}=\mathrm{O}_p(1) ,
\end{eqnarray*}
 and analogously for  $\tilde\sigma_{X,z}^2$.
\end{coro}

\begin{prop}\label{theorem2}
Under conditions \textup{[A0]--[A3]} in Section \ref{sec:estimation} and \textup{(i)--(ix)}, the
final estimates of $\sigma_X^2$ and $\sigma_Y^2$
(\ref{errorvarfinal}) converge in probability to their corresponding
true counterparts, that is,
\begin{eqnarray}
\hat\sigma_X^2\stackrel{P}{\rightarrow}\sigma_X^2 ,\qquad
\hat\sigma_Y^2\stackrel{P}{\rightarrow}
\sigma_Y^2.\nonumber
\end{eqnarray}
\end{prop}

\begin{pf} The result follows
straightforwardly from Corollary \ref{sigmacon}.
\end{pf}

While Lemma \ref{ymwjasalemma12dvar} implies consistency of the
estimator of the variance, we also require an extension regarding
estimation of the cross-covariance function.
Let $\{{\tilde\theta}_{p}(s, t, u, v), p=1, 2, \ldots, l\}$ be a
collection of real functions ${\tilde\theta}_p\dvtx
\mathbb{R}^4\rightarrow\mathbb{R}$.

\begin{enumerate}[{[C2.1c]}]
\item[{[C2.1c]}] For $\ell\geq|\nu|+2$ and any pair of $\ell_1$
and $\ell_2$ such
that $\ell=\ell_1+\ell_2$, $\ell_1\geq\nu_1+1$ and $\ell_2\geq
\nu_2+1$, we have, uniformly in $z\in\mathcal{Z}$, bandwidth $h_{1,z}$ and
$h_{2,z}$ satisfy $h_{1,z}\rightarrow0$,
$h_{1,z}/h_{2,z}\rightarrow1$, $n_{z,
h}h_{1,z}^{|\nu|+2}\rightarrow\infty$,
$n_{z,h}h_{1,z}^{2\ell+4}<\infty$ as $n\rightarrow\infty$.
\end{enumerate}
Define $\varrho_{p{\tilde\theta}, z}=\varrho_{p{\tilde\theta}, z}(s,
t)=\frac{\partial^{|\nu|}}{\partial s^{\nu_1}\, \partial t^{\nu_2}}
\int\int{\tilde\theta}_p(s, t, u, v)g_{XY,z}(s, t, u, v)\,\mathrm{d}u
\,\mathrm{d}v$ and
\begin{eqnarray*}
{\tilde\Theta}_{pn,z}&=&{\tilde\Theta}_{pn,z}(s, t)
\\
&=&\frac{1}{n_{z,h}
h_{1,z}^{\nu_1+1}h_{2,z}^{\nu_2+1}}
 \sum_{i\in\mathcal{N}_{z,h}}\frac{1}{EN} \sum_{1\leq j\leq
N_i}{\tilde\theta}_p(S_{ij}, T_{ij}, U_{ij}, V_{ij})
\kappa_2\biggl(\frac{S_{ij}-s}{h_{1,z}},
\frac{T_{ij}-t}{h_{2,z}}\biggr).
\end{eqnarray*}

\begin{lem}
Under conditions \textup{[A0]--[A3]}, \textup{(i),
(ii), (iii),
(viii), [C1.1b]} with $\mathcal{R}_1=\mathcal{S}$ and
$\mathcal{R}_2=\mathcal{T}$, \textup{[C1.2b]} with $g(\cdot, \cdot, \cdot,
\cdot)= g_{XY,z}(\cdot, \cdot, \cdot, \cdot)$ and \textup{[C2.1c]} (with
$\ell_1=\ell_2=1$ and $\nu_1=\nu_2=0$), we have $
\tilde\vartheta_{pn}=\sup_{(z, s, t)\in
\mathcal{Z}\times\mathcal{S}\times\mathcal{T}}
{\vert{\tilde\Theta}_{pn,z}(s,t)-\varrho_{p{\tilde\theta
},z}(s,t)\vert}/({h+(\sqrt{n_{z,h}}
h_{Y,1}^{\nu_1+1}h_{Y,2}^{\nu_2+1})^{-1}})=\mathrm{O}_p(1).$
\end{lem}

\begin{pf}
The proof is analogous to that of Lemmas 1 and
\ref{ymwjasalemma12dvar}.
\end{pf}

\begin{enumerate}[{[C6]}]
\item[{[C6]}] Uniformly in $z\in\mathcal{Z}$, bandwidths $h_{1,z}$
and $h_{2,z}$ satisfy $h_{1,z}\rightarrow0$,
$h_{1,z}/h_{2,z}\rightarrow1$,
$n_{z, h}h_{1,z}^{6}\rightarrow\infty$, $n_{z,h}h_{1,z}^{8}<\infty$
as $n\rightarrow\infty$.
\end{enumerate}

\begin{lem}[(Convergence of the cross-covariance function between $X$ and
$Y$)]\label{lemmacrosscov} Under conditions \textup{[A0]--[A3]}, \textup{(i),
(ii), (iii),
(viii), (ix),
[C3]} and
\textup{[C6]},
\[
\sup_{(z, s, t)\in\mathcal{Z}\times\mathcal{S}\times
\mathcal{T}}{\vert\tilde C_{XY, z}(s, t)- C_{XY, z}(s,
t)\vert}/{\bigl(h+\bigl(\sqrt{n_{z,h}}
h_{1,z}h_{2,z}\bigr)^{-1}\bigr)}=\mathrm{O}_p(1).
\]
\end{lem}

\begin{pf}
The proof is similar to that of Lemma \ref{ymwjasathm12dvar}.\vspace{2pt}
\end{pf}

Consider the real separable Hilbert space $L^2_Y(\mathcal{T})\equiv
H_Y$ (resp.~$L^2_X(\mathcal{S})\equiv H_X$) endowed with inner
product $\langle f, g\rangle_{H_Y}=\int_{\mathcal{T}}f(t)g(t)\,\mathrm{d}t$
(resp.~$\langle f, g\rangle_{H_X}=\int_{\mathcal{S}}f(s)g(s)\,\mathrm{d}s$) and
norm $\| f\|_{H_X}=\sqrt{\langle f, f\rangle_{H_X}}$
(resp.~$\| f\|_{H_Y}=\sqrt{\langle f, f\rangle_{H_Y}}$)
(\citet{couranthilbert1953}). Let $\mathcal{I}'_{Y,z}$ (resp.~$\mathcal
{I}'_{X,z}$) be the set of indices of the eigenfunctions
$\phi_{z,k}(t)$ (resp.~$\psi_{z,m}(s)$) corresponding to eigenvalues
$\lambda_{z,k}$ (resp.~$\rho_{z,m}$) of multiplicity one. We obtain
the consistency of $\tilde\lambda_{z,k}$ (resp.~$\tilde\rho_{z,m}$)
for $\lambda_{z,k}$ (resp.~$\rho_{z,m}$), the consistency of
$\tilde\phi_{z,k}(t)$ (resp.~$\tilde\psi_{z,m}(s)$) for
$\phi_{z,k}(t)$ (resp.~$\psi_{z,m}(s)$) in the $L^2_Y$-
(resp.~$L_X^2$-) norm $\|\cdot\|_{H_X}$
(resp.~$\|\cdot\|_{H_Y}$) when $\lambda_{z,k}$
(resp.~$\rho_{z,m}$) is of multiplicity one, and the uniform
consistency of
$\tilde\phi_{z,k}(t)$ (resp.~$\tilde\psi_{z,m}(s)$) for
$\phi_{z,k}(t)$ (resp.~$\psi_{z,m}(s)$) as well.

For $f, g, h\in H_Y$, define the rank one operator $f\otimes g\dvtx
h\rightarrow\langle f, h\rangle g$. Denote the separable Hilbert
space of Hilbert--Schmidt operators on $H_Y$ by $F_Y\equiv
\sigma_2(H_Y)$, endowed with $\langle T_1,
T_2\rangle_{F_Y}=\operatorname{tr}(T_1T_2^*)=\sum_{j}\langle T_1u_j, T_2
u_j\rangle_{H_Y}$ and $\| T\|_{F_Y}^2=\langle T,
T\rangle_{F_Y}$, where $T_1$, $T_2$, $T\in F_Y$, $T^*_2$ is
the\vspace{1pt}
adjoint of $T_2$ and $\{u_j\dvt j\geq1\}$ is any complete
orthonormal system in $H_Y$. The covariance operator $\mathbf{G}_{Y,z}$
(resp.~$\tilde\mathbf{G}_{Y,z}$) is generated by the kernel $G_{Y,z}$
(resp.~$\tilde G_{Y,z}$), that is,
$\mathbf{G}_{Y,z}(f)=\int_{\mathcal{T}}G_{Y,z}(t_1, t) f(t_1)\,\mathrm{d}t_1$
(resp.~$\tilde\mathbf{G}_{Y,z}(f)=\int_{\mathcal{T}}\tilde
G_{Y,z}(t_1, t)
f(t_1)\,\mathrm{d}t_1$). Obviously, $\mathbf{G}_{Y,z}$ and $\tilde\mathbf
{G}_{Y,z}$ are
Hilbert--Schmidt operators. As a result of (\ref{covyunifconrate}),
we have $\sup_{z\in\mathcal{Z}} \|\tilde\mathbf{G}_{Y,z}-
\mathbf{G}_{Y,z}\|_{F_{Y}}/(h+(\sqrt{n_{z,h}}
h_{Y,z}^2)^{-1})=\mathrm{O}_p(1)$.

Let $\mathcal{I}_{Y,z,i}=\{j\dvt \lambda_{z,j}=\lambda_{z,i}\}$ and
$\mathcal{I}_{Y,z}'=\{i\dvt |\mathcal{I}_{Y,z,i}|=1\}$,
where $|\mathcal{I}_{Y,z,i}|$ denotes the number of
elements in $\mathcal{I}_{Y,z,i}$. Define
$\mathbf{P}_{z,j}^Y=\sum_{k\in\mathcal{I}_{Y,z,j}}\phi
_{z,k}\otimes
\phi_{z,k}$ and $\tilde
\mathbf{P}_{z,j}^Y=\sum_{k\in\mathcal{I}_{Y,z,j}}\tilde\phi
_{z,k}\otimes
\tilde\phi_{z,k}$ to be the true and estimated orthogonal projection
operators from $H_Y$ to the subspace spanned by $\{\phi_{z,k}\dvt k\in
\mathcal{I}_{Y,z,j}\}$. Set $\delta_{z,
j}^{Y}=\frac{1}{2}\min\{|\lambda_{z,l}-\lambda_{z,j}|\dvt
l\notin\mathcal{I}_{Y,z,j}\}$
and $\bolds\Lambda_{\delta_{z, j}^{Y}}=\{c\in\mathcal{C}\dvt
|c-\lambda_{z,j}|=\delta_{z, j}^{Y}\}$, where
$\mathcal{C}$ stands for the complex numbers. Let $\mathbf{R}_{Y,z}$
(resp.~$\tilde\mathbf{R}_{Y,z}$) be the resolvent of $\mathbf{G}_{Y,z}$
(resp.~$\tilde\mathbf{G}_{Y,z}$), that is, $\mathbf
{R}_{Y,z}(c)=(\mathbf{G}_{Y,z}-c
I)^{-1}$ (resp.~$\tilde\mathbf{R}_{Y,z}(c)=(\tilde\mathbf
{G}_{Y,z}-c I)^{-1}$). Let
$A_{\delta_{z, j}^{Y}}=\sup\{\|
\mathbf{R}_{Y,z}(c)\|_{F_{Y}}\dvt c\in\bolds\Lambda
_{\delta_{z,
j}^{Y}}\}$ and\vspace{2pt}
%
\begin{equation}\label{definofalpha}
\alpha_X ={( \delta_{z,j}^X (A_{\delta_{z,j}^X})^2
)}/{\bigl(\bigl(h+\bigl(\sqrt{n_{z,h}}h_{X,z}^2\bigr)^{-1}\bigr)^{-1}-A_{\delta
_{z,j}^X}\bigr)}.\vspace{2pt}
\end{equation}
Parallel notation is assumed for the $Y$ process.\vspace{2pt}

\begin{prop} \label{theorem1}
Under conditions \textup{[A0]--[A3]} in Section \ref{sec:estimation} and
conditions \textup{(i)--(iii),
(viii), (ix), [C3], [C4]} and \textup{[C6]},
it holds that\vspace{2pt}
%
\begin{eqnarray}\label{xeigenvalueconv}
|\tilde\rho_{z,m} -\rho_{z,m}|&=&\mathrm{O}_p(\alpha_X),
\\[2pt]\label{xeigenfunconv}
\|\tilde\psi_{z,m}
-\psi_{z,m}\|_{H_X}&=&\mathrm{O}_p(\alpha_X) ,\qquad
m\in\mathcal{I}_{X, z}' ,
\\[2pt]\label{xeigenfununifconv}
\sup_{s\in\mathcal{S}}|\tilde\psi_{z,m}(s)
-\psi_{z,m}(s)|&=&\mathrm{O}_p(\alpha_X) ,    \qquad        m\in
\mathcal{I}_{X, z}' ,\vadjust{\goodbreak}
\end{eqnarray}
%
\begin{eqnarray}\label{yeigenvalueconv}
|\tilde\lambda_{z,k} -\lambda_{z,k}|&=&\mathrm{O}_p(\alpha_Y),
\\\label{yeigenfunconv}
\|\tilde\phi_{z,k}
-\phi_{z,k}\|_{H_Y}&=&\mathrm{O}_p(\alpha_Y) ,\qquad
k\in\mathcal{I}_{Y, z}' ,
\\\label{yeigenfununifconv}
\sup_{t\in\mathcal{T}}|\tilde\phi_{z,k}(t)
-\phi_{z,k}(t)|&=&\mathrm{O}_p(\alpha_Y) ,   \qquad         k\in
\mathcal{I}_{Y, z}' ,
\end{eqnarray}
%
\begin{eqnarray}\label{xyeigenvalueconv}
|\tilde\sigma_{z,mk}-\sigma_{z, mk}|&=&\mathrm{O}_p\bigl(\max\bigl(\alpha_X,
\alpha_Y, h+\bigl(\sqrt{n_{z,h}} h_{1,z}h_{2,z}\bigr)^{-1}\bigr)\bigr) ,
\end{eqnarray}
where the norms on $H_X$ and $H_Y$ are defined on page 29,
both $\alpha_X, \alpha_Y$ are defined in \textup{(\ref{definofalpha})}
and converge to zero as $n\rightarrow\infty$ and the above $O_p$
terms are uniform in $z\in\mathcal{Z}$.
\end{prop}

\begin{pf} The proof is similar to
the proof of Theorem 2 in \citet{YMW2005jasa}. The uniformity
result follows from that of Lemmas \ref{ymwjasathm12dvar} and
\ref{lemmacrosscov}.
\end{pf}

Note that
%
\begin{eqnarray}\label{betadefn}
\beta(z, s,t)=\sum_{k=1}^\infty\sum_{m=1}^\infty
\frac{E(\zeta_{z,m}\xi_{z,k})}{E(\zeta_{z,m}^2)}\psi_{z,m}(s)\phi
_{z,k}(t).
\end{eqnarray}
To define the convergence of the right-hand side of (\ref{betadefn}),
in the $L_2$ sense, in $(s,t)$ and uniformly in $z$, we require that
[A4] $\sum_{k=1}^\infty\sum_{m=1}^\infty
\sigma_{z,mk}^2/\rho_{z,m}^2 < \infty$ uniformly for
$z\in\mathcal{Z}$.

The proof of the following result is straightforward.
\begin{lem}Under condition
\textup{[A4]}, uniformly in $z\in\mathcal{Z}$, the right-hand side of
\textup{(\ref{betadefn})} converges in the $L_2$ sense.
\end{lem}

The next
result is stated without proof and requires assumptions [A4] and the
following:
\begin{eqnarray*}
\sum_{m=1}^{M(n)}\frac{\delta_{z,m}^X
(A_{\delta_{z,m}^X})^2
}{(h+(\sqrt{n_{z,h}}h_{X,z}^2)^{-1})^{-1}-A_{\delta
_{z,m}^X}}&\rightarrow&0 , \nonumber
\\
\hspace{-26pt}\mbox{[A5]}  \qquad\sum_{k=1}^{K(n)} \frac{ \delta_{z,k}^Y
(A_{\delta_{z,k}^Y})^2
}{(h+(\sqrt{n_{z,h}}h_{Y,z}^2)^{-1})^{-1}-A_{\delta
_{z,k}^Y}}&\rightarrow
& 0 \qquad        \mbox{uniformly in $z\in\mathcal{Z}$,}
\\
MK\bigl(h+\bigl(\sqrt{n_{z,h}} h_{1,z}h_{2,z}\bigr)^{-1}\bigr)&\rightarrow& 0 .
\end{eqnarray*}

\begin{lem}\label{rawbetaconsistent}
Under conditions of Proposition \ref{theorem1}, \textup{[A4]} and \textup{[A5]},
%
\begin{eqnarray}\label{rawbetaconvinl2}
\lim_{n\rightarrow\infty}\sup_{z\in\mathcal{Z}}
\int_{\mathcal{S}}\int_{\mathcal{T}} [\tilde\beta(z, s,t)-\beta(z,
s, t)]^2=0    \qquad       \mbox{in
probability}.
\end{eqnarray}
\end{lem}

\begin{pf*}{Proof of Theorem \ref{theorem:finalbetaconv}} We consider
only the convergence of $\hat\beta(z,s,t)$. The consistency of
$\hat\mu_{X,z}(s)$ and $\hat\mu_{Y,z}(t)$ is analogous. First, note
that
%
\begin{eqnarray}\label{betarateupbd}
&&\int_{\mathcal{T}}\int_{\mathcal{S}} \bigl(\hat\beta(z, s, t)-\beta
(z, s, t)\bigr)^2\,\mathrm{d}s\,\mathrm{d}t\nonumber
\\
&&\quad \leq2(2b/h+1)\sum_{p=1}^{P}\omega_{0,2}\bigl(z^{(p)}, z, b\bigr)^2
\int_{\mathcal{T}}\int_{\mathcal{S}} \bigl(\tilde\beta\bigl( z^{(p)}, s,
t\bigr)-\beta\bigl(z^{(p)},s,t\bigr)\bigr)^2\,\mathrm{d}s\,\mathrm{d}t\qquad
\\
&&{}\qquad +2\int_{\mathcal{T}}\int_{\mathcal{S}}
\Biggl(\sum_{p=1}^{P}\omega_{0,2}\bigl(z^{(p)}, z, b\bigr) \beta\bigl(z^{(p)},s,t\bigr)-\beta
(z, s,
t)\Biggr)^2\,\mathrm{d}s\,\mathrm{d}t,\nonumber
\end{eqnarray}
where the $2b/h+1$ in the last inequality is due to the fact that the
kernel function $K(\cdot)$ is of bounded support $[-1, 1]$.
Let $a(k)=\sum_{p=1}^P K_b(z^{(p)}-z) (z^{(p)}-z)^k$,
$b(k)=\sum_{p=1}^P K_b(z^{(p)}-z)^2 (z^{(p)}-z)^k$, $\mu_k=\int
K(u)u^k\,\mathrm{d}u$ and $\nu_k=\int(K(u))^2 u^k\,\mathrm{d}u$. We then have
\[
a(k)=\mu_k \frac{b^k}{h}\bigl(1+\mathrm{o}(1)\bigr) \quad\mbox{and}\quad  b(k)=\nu_k
\frac{b^{k-1}}{h}\bigl(1+\mathrm{o}(1)\bigr)
\]
for small $h$ (large $P\propto1/h$) and small $b$. Moreover, the usual
boundary techniques can be applied near the two end points.
Consequently,we have
\begin{eqnarray*}
\sum_{p=1}^{P} \omega_{0,2}\bigl(z^{(p)}, z, b\bigr)^2&=&\mathbf{e}_{1,2}^T
(\mathbf{C}
^T\mathbf{W}\mathbf{C})^{-1}
(\mathbf{C}^T\mathbf{W}\mathbf{W}\mathbf{C}) (\mathbf{C}^T\mathbf
{W}\mathbf{C})^{-1} \mathbf{e}_{1,2}\nonumber
\\
&=& \mathbf{e}_{1,2}^T \left(
\matrix{
a(0) & a(1)\cr
a(1) & a(2)
}
\right)^{-1}
\left(
\matrix{
b(0) & b(1)\cr
b(1) & b(2)
}
\right)
\left(
\matrix{
a(0) & a(1)\cr
a(1) & a(2)
}
\right)^{-1}
\mathbf{e}_{1,2}
\\
&=&\biggl(\frac{\mu_2^2\nu_0-2\mu_1\mu_2\nu_1+\mu_1^2\nu_2}{\mu_0\mu
_2-\mu_1^2}\biggr)
\biggl(\frac{b}{h}\biggr)\bigl(1+\mathrm{o}(1)\bigr).\nonumber
\end{eqnarray*}
Due to the compactness of $\mathcal{Z}$, the above o-term is uniform
in $z\in\mathcal{Z}$. This implies that
%
\begin{eqnarray}\label{weightsqintconv}
\int_{\mathcal{Z}} \sum_{p=1}^{P} \omega_{0,2}\bigl(z^{(p)}, z, b\bigr)^2\, \mathrm{d}z=
\biggl(\frac{\mu_2^2\nu_0-2\mu_1\mu_2\nu_1+\mu_1^2\nu_2}{\mu_0\mu
_2-\mu_1^2}\biggr)
\biggl(\frac{b}{h}\biggr)|\mathcal{Z}|\bigl(1+\mathrm{o}(1)\bigr)
\end{eqnarray}
for small $h$ and $b$, where $|\mathcal{Z}|$ denotes the
Lebesgue measure of $\mathcal{Z}$. Hence, (\ref{weightsqintconv})
and the consistency of $\tilde\beta(z,s,t)$ in the $L_2$ sense in
$(s,t)$ and uniformly in $z$ due to (\ref{rawbetaconvinl2})
imply that
%
\begin{eqnarray}\label{betaconv1part}
\int_{\mathcal{Z}}\Biggl[\sum_{p=1}^{P}\omega_{0,2}\bigl(z^{(p)}, z, b\bigr)^2
\int_{\mathcal{T}}\int_{\mathcal{S}} \bigl( \bigl(\tilde\beta\bigl( z^{(p)}, s,
t\bigr)-\beta\bigl(z^{(p)},s,t\bigr)\bigr)\bigr)^2\,\mathrm{d}s\,\mathrm{d}t\Biggr]\,\mathrm{d}z\stackrel{P}{\rightarrow}0.\quad
\end{eqnarray}
For the second part in (\ref{betarateupbd}), applying a Taylor
expansion of $\beta(z^{(p)},s,t)$ at each $z$,
we have
\begin{eqnarray*}
&&\sum_{p=1}^{P}\omega_{0,2}\bigl(z^{(p)}, z, b\bigr) \beta
\bigl(z^{(p)},s,t\bigr)\nonumber
\\
&&\quad =\mathbf{e}_{1,2}^T \left(
\matrix{
a(0) & a(1)\cr
a(1) & a(2)
}
\right)^{-1}\left(
\matrix{
a(0)\cr
a(1)
}
\right) \beta(z, s, t)+
\mathbf{e}_{1,2}^T \left(
\matrix{
a(0) & a(1)\cr
a(1) & a(2)
}
\right)^{-1}\left(
\matrix{
a(1)\cr
a(2)
}
\right) \frac{\partial}{\partial z}\beta(z, s,
t)\nonumber
\\
&&{}\qquad +\frac{1}{2}\mathbf{e}_{1,2}^T \left(
\matrix{
a(0) & a(1)\cr
a(1) & a(2)
}
\right)^{-1}\left(
\matrix{
a(2)\cr
a(3)
}
\right) \frac{\partial^2}{\partial z^2}\beta(z, s, t) +
{\mbox{ higher order terms}}\nonumber
\\
&&\quad =\beta(z,s,t)+\frac{1}{2}b^2\frac{\mu_2^2-\mu_1\mu_3}{\mu_0\mu
_2-\mu_1^2}\frac{\partial^2}{\partial
z^2}\beta(z, s, t)+{\mbox{ higher order terms}}.
\end{eqnarray*}
Hence, $\sum_{p=1}^{P}\omega_{0,2}(z^{(p)}, z, b)
\beta(z^{(p)},s,t)-\beta(z, s,
t)=\frac{1}{2}b^2\frac{\mu_2^2-\mu_1\mu_3}{\mu_0\mu_2-\mu
_1^2}\frac{\partial^2}{\partial
z^2}\beta(z, s, t)(1+\mathrm{o}(1))$ and
\begin{eqnarray}\label{betaconv2part}
&&\int_{\mathcal{Z}}\int_{\mathcal{T}}\int_{\mathcal{S}}
\Biggl(\sum_{p=1}^{P}\omega_{0,2}\bigl(z^{(p)}, z, b\bigr) \beta\bigl(z^{(p)},s,t\bigr)-\beta
(z, s,
t)\Biggr)^2\,\mathrm{d}s\,\mathrm{d}t\,\mathrm{d}z\nonumber
\\[-8pt]\\[-8pt]
&&\quad =\frac{1}{2}b^2\frac{\mu_2^2-\mu_1\mu_3}{\mu_0\mu_2-\mu
_1^2}\biggl(\int_{\mathcal{Z}}\int_{\mathcal{T}}\int_{\mathcal{S}}
\frac{\partial^2}{\partial z^2}\beta(z, s,
t)\,\mathrm{d}s\,\mathrm{d}t\,\mathrm{d}z\biggr)\bigl(1+\mathrm{o}(1)\bigr)\rightarrow0.\nonumber
\end{eqnarray}
Combining (\ref{betaconv1part}) and (\ref{betaconv2part}), and
further noting condition (xi), completes
the proof.
\end{pf*}

\begin{pf*}{Proof of Theorem \ref{theorem4}} Note that
\begin{eqnarray*}
Y^*(t)-\hat
Y^*(t)&=&\mu_{Y,Z^*}(t)-\hat\mu_{Y,Z^*}(t)+\int_{\mathcal
{S}}\bigl(\beta(Z^*,
s, t)-\hat\beta(Z^*, s, t)\bigr) \bigl(X^*(s)-\mu_{X,
Z^*}(s)\bigr)\,\mathrm{d}s
\\
&&{}-\int_{\mathcal{S}}\hat\beta(Z^*, s, t) \bigl(\mu_{X,
Z^*}(s)-\hat\mu_{X, Z^*}(s)\bigr)\,\mathrm{d}s.
\end{eqnarray*}
The convergence results in Theorem \ref{theorem:finalbetaconv}
imply that $\int_{\mathcal{T}} (Y^*(t)-\hat Y^*(t))^2\,\mathrm{d}t\stackrel
{P}{\rightarrow}0$,
as desired.
\end{pf*}

\section*{Acknowledgements}
We wish to thank two referees for helpful comments. Yichao Wu's
research has been supported in part by NIH Grant R01-GM07261
and
NFS Grant DMS-09-05561. Jianqing Fan's research has been supported in part by National
Science Foundation (NSF) Grants DMS-03-54223 and DMS-07-04337.
Hans-Georg M\"{u}ller's research has been supported in part by National Science
Foundation (NSF) Grants DMS-03-54223, DMS-05-05537 and
DMS-08-06199.

\printhistory


\begin{thebibliography}{29}

\bibitem[\protect\citeauthoryear{Besse and Ramsay}{1986}]{MR848110}
{Besse, P.} and {Ramsay, J.O.} (1986).
Principal components analysis of sampled functions.
\textit{Psychometrika} \textbf{51} 285--311.
\MR{0848110}

\bibitem[\protect\citeauthoryear{Cai and Hall}{2006}]{CH2006}
{Cai, T.} and {Hall, P.} (2006).
Prediction in functional linear regression.
\textit{Ann. Statist.} \textbf{34} 2159--2179.
\MR{2291496}

\bibitem[\protect\citeauthoryear{Cardot}{2007}]{Cardot2007}
{Cardot, H.} (2007).
Conditional functional principal components analysis.
\textit{Scand. J. Statist.} \textbf{34} 317--335.
\MR{2346642}

\bibitem[\protect\citeauthoryear{Cardot \textit{et al.}}{2003}]{CFS2003SS}
{Cardot, H.}, {Ferraty, F.} and {Sarda, P.} (2003).
Spline estimators for the functional linear model.
\textit{Statist. Sin.} \textbf{13} 571--591.
\MR{1997162}

\bibitem[\protect\citeauthoryear{Cardot and Sarda}{2008}]{cardotsarda2005}
{Cardot, H.} and {Sarda, P.} (2008).
Varying-coefficient functional linear regression models.
\textit{Comm. Statist. Theory Methods} \textbf{37} 3186--3203.

\bibitem[\protect\citeauthoryear{Carey \textit{et~al.}}{1998}]{Carey1998}
{Carey, J.R.}, {Liedo, P.}, {M\"{u}ller, H.G.}, {Wang,
J.L.} and {Chiou, J.M.} (1998).
Relationship of age patterns of fecundity to mortality, longevity,
and lifetime reproduction in a large cohort of mediterranean fruit fly
females.
\textit{Journals of Gerontology Series A: Biological
Sciences and Medical Sciences} \textbf{53}
245--251.

\bibitem[\protect\citeauthoryear{Courant and
Hilbert}{1953}]{couranthilbert1953}
{Courant, R.} and {Hilbert, D.} (1953).
\textit{Methods of Mathematical Physics}. New York: Wiley.

\bibitem[\protect\citeauthoryear{Cuevas \textit{et al.}}{2002}]{CFF2002}
{Cuevas, A.}, {Febrero, M.} and {Fraiman, R.} (2002).
Linear functional regression: The case of fixed design and functional
response.
\textit{Canad. J. Statist.} \textbf{30} 285--300.
\MR{1926066}

\bibitem[\protect\citeauthoryear{Fan and Gijbels}{1996}]{FG1996}
{Fan, J.} and {Gijbels, I.} (1996).
\textit{Local Polynomial Modelling and Its Applications}. London:
Chapman \& Hall.
\MR{1383587}

\bibitem[\protect\citeauthoryear{Fan \textit{et al.}}{2007}]{FHL2007}
{Fan, J.}, {Huang, T.} and {Li, R.} (2007).
Analysis of longitudinal data with semiparametric estimation of
covariance function.
\textit{J. Amer. Statist. Assoc.} \textbf{35}
632--641.
\MR{2370857}

\bibitem[\protect\citeauthoryear{Fan and Zhang}{2000}]{FZ2000}
{Fan, J.} and {Zhang, J.-T.} (2000).
Two-step estimation of functional linear models with applications to
longitudinal data.
\textit{J. Roy. Statist. Soc. Ser. B} \textbf{62}
303--322.
\MR{1749541}

\bibitem[\protect\citeauthoryear{Hall and Horowitz}{2007}]{HH2007}
{Hall, P.} and {Horowitz, J.L.} (2007).
Methodology and convergence rates for functional linear regression.
\textit{Ann. Statist.} \textbf{35} 70--91.
\MR{2332269}

\bibitem[\protect\citeauthoryear{Hall \textit{et al.}}{2006}]{HallMuellerWang2006}
{Hall, P.}, {M\"{u}ller, H.} and {Wang, J.} (2006).
Properties of principal component methods for functional and
longitudinal data analysis.
\textit{Ann. Statist.} \textbf{34} 1493--1517.
\MR{2278365}

\bibitem[\protect\citeauthoryear{James \textit{et al.}}{2000}]{MR1789811}
{James, G.M.}, {Hastie, T.J.} and {Sugar, C.A.} (2000).
Principal component models for sparse functional data.
\textit{Biometrika} \textbf{87} 587--602.
\MR{1789811}

\bibitem[\protect\citeauthoryear{Liu and M\"{u}ller}{2009}]{Liumueller2009}
{Liu, B.} and {M\"{u}ller, H.} (2009).
Estimating derivatives for samples of sparsely observed functions,
with application to on-line auction dynamics.
\textit{J. Amer. Statist. Assoc.}
\textbf{104} 704--717.

\bibitem[\protect\citeauthoryear{M\"{u}ller \textit{et~al.}}{2001}]{MCWLV2001}
{M\"{u}ller, H.G.}, {Carey, J.R.}, {Wu, D.}, {Liedo,
P.} and {Vaupel, J.W.} (2001).
Reproductive potential predicts longevity of female Mediterranean
fruit flies.
\textit{Proc. Roy. Soc. Ser. B} \textbf{268} 445--450.

\bibitem[\protect\citeauthoryear{Pearson \textit{et~al.}}{1997}]{PMBLF1997}
{Pearson, J.D.}, {Morrell, C.H.}, {Brant, L.J.},
{Landis, P.K.} and {Fleg, J.L.} (1997).
Age-associated changes in blood pressure in a longitudinal study of
healthy men and women. \textit{Journals of Gerontology Series A: Biological
Sciences and Medical Sciences}, \textbf{52}
177--183.

\bibitem[\protect\citeauthoryear{Ramsay and
Dalzell}{1991}]{ramsaydalzell1991}
{Ramsay, J.} and {Dalzell, C.J.} (1991).
Some tools for functional data analysis (with discussion).
\textit{J. Roy. Statist. Soc. Ser. B} \textbf{53}
539--572.
\MR{1125714}

\bibitem[\protect\citeauthoryear{Ramsay and Silverman}{2002}]{RS2002}
{Ramsay, J.O.} and {Silverman, B.W.} (2002).
\textit{Applied Functional Data Analysis: Methods and Case Studies}.
New York: Springer.
\MR{1910407}

\bibitem[\protect\citeauthoryear{Ramsay and Silverman}{2005}]{RS2005}
{Ramsay, J.O.} and {Silverman, B.W.} (2005).
\textit{Functional Data Analysis}. New York:
Springer.
\MR{2168993}

\bibitem[\protect\citeauthoryear{Rice and Wu}{2001}]{ricewu2001}
{Rice, J.} and {Wu, C.} (2001).
Nonparametric mixed effects models for unequally sampled noisy
curves.
\textit{Biometrics} \textbf{57} 253--259.
\MR{1833314}

\bibitem[\protect\citeauthoryear{Rice and Silverman}{1991}]{MR1094283}
{Rice, J.A.} and {Silverman, B.W.} (1991).
Estimating the mean and covariance structure nonparametrically when
the data are curves.
\textit{J. Roy. Statist. Soc. Ser. B} \textbf{53}
233--243.
\MR{1094283}

\bibitem[\protect\citeauthoryear{Shock \textit{et~al.}}{1984}]{Shock1984}
{Shock, N.W.}, {Greulich, R.C.}, {Andres, R.},
{Lakatta, E.G.}, {Arenberg, D.} and {Tobin, J.D.}
(1984).
Normal human aging: The Baltimore longitudinal study of aging.
NIH Publication No. 84-2450,  U.S. Government
Printing Office., Washington, DC.

\bibitem[\protect\citeauthoryear{Silverman}{1996}]{MR1389877}
{Silverman, B.W.} (1996).
Smoothed functional principal components analysis by choice of norm.
\textit{Ann. Statist.} \textbf{24} 1--24.
\MR{1389877}

\bibitem[\protect\citeauthoryear{Staniswalis and
Lee}{1998}]{staniswalislee1998}
{Staniswalis, J.-G.} and {Lee, J.-J.} (1998).
Nonparametric regression analysis of longitudinal data.
\textit{J. Amer. Statist. Assoc.} \textbf{93}
1403--1418.
\MR{1666636}

\bibitem[\protect\citeauthoryear{Wahba}{1990}]{wahba1990}
{Wahba, G.} (1990).
\textit{Spline Models for Observational Data}. Philadelphia, PA:
Society for Industrial and Applied Mathematics.
\MR{1045442}

\bibitem[\protect\citeauthoryear{Yao \textit{et~al.}}{2005a}]{YMW2005jasa}
{Yao, F.}, {M\"{u}ller, H.-G.} and {Wang, J.-L.} (2005a).
Functional data analysis for sparse longitudinal data.
\textit{J.~Amer. Statist. Assoc.}
\textbf{100} 577--590.
\MR{2160561}

\bibitem[\protect\citeauthoryear{Yao \textit{et al.}}{2005b}]{YMW2005annal}
{Yao, F.}, {M\"{u}ller, H.-G.} and {Wang, J.-L.} (2005b).
Functional linear regression analysis for longitudinal data.
\textit{Ann. Statist.} \textbf{33} 2873--2903.
\MR{2253106}

\bibitem[\protect\citeauthoryear{Zhang \textit{et~al.}}{2008}]{ZLJOD2007}
{Zhang, C.M.}, {Lu, Y.F.}, {Johnstone, T.}, {Oaks, T.}
and {Davidson, R.J.} (2008).
Efficient modeling and inference for event-related functional fMRI
data.
\textit{Comput. Statist. Data. Anal.} \textbf{52} 4859--4871.

\end{thebibliography}
\end{document}